\documentclass[reqno]{amsart}

\usepackage{amsthm}
\usepackage{amssymb}
\usepackage{amsmath}

\newtheorem{theorem}{Theorem}[section]
\newtheorem*{theorem*}{Theorem}
\newtheorem{lemma}{Lemma}

\newtheorem*{proposition*}{Proposition}

\newtheorem{corollary}{Corollary}[section]
\newtheorem*{corollary*}{Corollary}
\theoremstyle{definition}
\newtheorem*{definition}{Definition}
\newtheorem*{comments*}{Comments}
\newtheorem{example}{Example}
\newtheorem*{example*}{Example}

\newtheorem*{remark*}{Remark}
\newtheorem*{remarks*}{Remarks}

\numberwithin{equation}{section}

\gdef\myletter{}

\let\savetheequation\theequation
\def\theequation{\savetheequation\myletter}

\newcommand{\CC}{{\mathbb C}}
\newcommand{\LL}{{\mathbb L}}
\newcommand{\PP}{{\mathbb P}}
\newcommand{\RR}{{\mathbb R}}

\newcommand{\fg}{{\mathfrak g}}
\newcommand{\g}{{\mathfrak g}}
\newcommand{\fh}{{\mathfrak h}}

\newcommand{\p}{{\mathfrak p}}

\newcommand{\fk}{{\mathfrak k}}

\newcommand{\K}{{\mathcal K}}

\renewcommand{\P}{{\mathcal P}}

\newcommand{\Identity}{\mbox{Identity}}
\newcommand{\Res}{\mbox{Res}}

\def    \to     {{\longrightarrow}}

\begin{document}

\title[Equivariant de Rham theory and graphs]
{Equivariant de Rham theory and graphs}
\author[V. Guillemin]{V. Guillemin\footnotemark {$^\dagger$}}
\thanks{Supported by NSF grant DMS 890771}
\address{Department of Mathematics, MIT, Cambridge, MA 02139}
\email{vwg@@math.mit.edu}
\thanks{printed \today}
\author[C. Zara]{C. Zara}
\address{Department of Mathematics, MIT, Cambridge, MA 02139}
\email{czara@math.mit.edu}

\begin{abstract}
Goresky, Kottwitz and MacPherson have recently shown that 
the computation of the equivariant cohomology ring of a  
$G$-manifold can be reduced to a computation in graph 
theory. This opens up the possibility that many of the 
fundamental theorems in equivariant de Rham theory may, 
on closer inspection, turn out simply to be theorems about 
graphs. In this paper we show that for some familiar theorems,  
this is indeed the case. 
\end{abstract}

\maketitle

\subsection*{Introduction}
\label{intro}

This article will consist of two essentially disjoint parts.
Part~1 is an exposition of (mostly) well-known results about
$G$-manifolds.  In section~\ref{sec:1.1}--\ref{sec:1.3} we review the
  definition of the equivariant de Rham cohomology 
ring of a $G$-manifold and recall the statements of the 
two fundamental ``localization theorems'' in equivariant
de Rham theory: the Atiyah-Bott-Berline-Vergne theorem and the 
Jeffrey-Kirwan theorem. 
In section~\ref{sec:1.4} we discuss the
  ``Smith'' problem for $G$-manifolds (which is
  concerned with the question:  Given a $G$-manifold 
with isolated fixed points, what kinds of
representations can occur as isotropy representations at the
fixed points?)  Then in sections~\ref{sec:1.5}--\ref{sec:1.7} we 
report on some very exciting recent results of 
Goresky-Kottwitz-MacPherson which have to do with the tie-in between
``equivariant de~Rham theory'' and ``graphs'' alluded to in our
title. These results show that for a large class of $G$-manifolds, $M$, 
with $M^G$ finite, the equivariant cohomology ring of $M$ is 
isomorphic to 
the equivariant cohomology ring of a pair
$(\Gamma, \alpha)$, where 
$\Gamma$ is the intersection graph of a necklace of embedded
$S^2$'s, each of which is equipped with a circle action (i.e.,~an 
axis of symmetry), and $\alpha$ is an ``axial'' function which
describes the directions in which the axes of these $S^2$'s are tilted. Finally, in section \ref{sec:1.8} we discuss a Morse theoretic
recipe for computing the Betti numbers of $M$ in terms of the pair
$(\Gamma, \alpha)$.

The second part of this article is concerned with the
combinatorial invariants of a pair $(\Gamma, \alpha)$, $\Gamma$
being {\em any} finite simple $d$-valent graph and $\alpha$ an 
abstract analogue of the axial function alluded to above. 
In particular, for such a pair we will prove combinatorial
versions of the theorems described in sections 
\ref{sec:1.2}-\ref{sec:1.3} and \ref{sec:1.8}.
These combinatorial ``localization'' theorems help to shed some
light on the role of the localization theorems in 
Smith theory : 
From the localization theorems one can
generate a lot of complicated identities among the
weights of the isotropy representations.  However, the question
of whether one can extract from these identities any new
information about the isotropy representations themselves has
been an open question for a long time.  Our graph theoretical
results seem to indicate that one can't.

This article is the first of a series of two articles on graphs
and equivariant cohomology. In the second article in this series 
we will discuss K-theoretical analogues of the results above
and give a purely combinatorial proof of the so-called 
``quantization commutes with reduction'' conjecture.

\setcounter{section}{1}
\subsection{Equivariant de Rham theory}
\label{sec:1.1}

Let $G$ be an $n$-dimensional Lie group which is compact,
connected and abelian, i.e.,~an $n$-dimensional torus.  Let $\fg$ 
be its Lie algebra and $\fg^*$ the vector space dual of $\fg$.
We will fix a basis $\xi^1, \ldots , \xi^n$ of $\fg$ and let
$x_1, \ldots , x_n$ be the dual basis.  Using this basis, the
symmetric algebra $S(\fg^*)$ can be identified with the
polynomial ring $\CC [x_1, \ldots, x_n]$.  

Let $M$ be a $2d$-dimensional manifold and $\tau$ an action of
$G$ on $M$.  From $\tau $ one gets an infinitesimal action,
$\partial \tau$, of $\fg$ on $M$ which associates to every
element $\xi$ of $\fg$ a vector field $\xi_M$.  Let $\Omega
(M)$ be the usual complex of de~Rham forms on $M$ and $\Omega
(M)^G$ the subcomplex of $G$-invariant de~Rham forms.  One defines
the \emph{equivariant} de~Rham complex of 
$M$ to be the tensor product
\begin{equation}
  \label{eq:1.1}
  \Omega_G (M) = \Omega (M)^G \otimes S (\fg^*)
\end{equation}
with the coboundary operator
\begin{equation}
  \label{eq:1.2}
  d_G (\omega \otimes f) = d \omega \otimes f + \sum \iota
  (\xi^i_M) \omega \otimes x_i f \, .
\end{equation}
The \emph{equivariant cohomology ring} of $M, H_G (M)$, is the
cohomology ring of this complex.  A few properties of this ring
which we will need below are:

\begin{enumerate}
\item $H_G(M)$ is an $S(\fg^*)$-module.  (This follows from the
  fact that $\Omega_G(M)$ is an $S(\fg^*)$ module by
  (\ref{eq:1.1}) and $d_G$ is an $S(\fg^*)$ module morphism by
  (\ref{eq:1.2}).)

\item  $H_G (pt) = S (\fg^*)$.

\item  Suppose $M$ is compact and oriented.  Then there is an
  integration operation
\begin{equation}
  \label{eq:1.3}
  \int : \Omega_G (M) \to S(\fg^*)
\end{equation}
defined by 
\begin{displaymath}
  \int (\omega \otimes f) = f \int \omega \, .
\end{displaymath}
It is easily checked that $\int d_G=0$ and hence that this
integration operation induces an integration operation on
cohomology
\begin{equation}
\label{eq:1.4}
\int  : H_G (M) \to S(\fg^*)  \, .
\end{equation}

\item  One can write $d_G$ as a sum, $d_1 + d_2$, $d_1$ and $d_2$ 
  being the first and second terms on the right hand side of
  (\ref{eq:1.2}).  Thus $\Omega_G (M)$ is a \emph{bi-complex},
  and the additive structure of $H_G (M)$ can be computed by the
  spectral sequence of this bi-complex.  The $E_1$ term in this
  spectral sequence is the $d_1$-cohomology of $\Omega_G(M)$,
  namely
  \begin{equation}
    \label{eq:1.5}
    H(M) \otimes S (\fg^*) \, .
  \end{equation}
One says that $M$ is \emph{equivariantly formal} if the spectral
sequence is trivial, i.e.~if, as vector spaces,
\begin{equation}
  \label{eq:1.6}
  H_G(M)=H(M) \otimes S(\fg^*) \, .
\end{equation}
\end{enumerate}
One can show, by the way, that if (\ref{eq:1.6}) holds as an
identity of vector spaces, it also holds as an identity of
$S(\fg^*)$-modules.  However, (\ref{eq:1.6}) doesn't, in general, 
tell one very much about the \emph{ring} structure of $H_G(M)$
(about which we will have more to say in \S\ref{sec:1.6}).

The property of being equivariantly formal is a bit technical;
however there are a number of interesting assumptions on $M$ which 
will imply this property.  (See~\cite{GKM}.)  Of these
assumptions the one that will be of most interest to us is the
following:

\begin{theorem*}[Kirwan]
  \label{Kirwan}

If $M$ is a symplectic manifold and the action $\tau$ is
Hamiltonian, $M$ is equivariantly formal.

\end{theorem*}

\subsection{The Atiyah-Bott-Berline-Vergne localization theorem}
\label{sec:1.2}

Let $M$ be compact and oriented and, also, to simplify the
statement of the localization theorem, let $M^G$ be finite.  For
$p \in M^G$ one has an isotropy representation $\tau_p$ of $G$
on $T_p$ and we will denote the weights of this representation 
by $\alpha_{i,p}, i=1 , \ldots , d$.  Since $\tau_p$ is a
\emph{real} representation, these weights are, strictly speaking, 
only defined up to sign; however, since $M$ is oriented, the
product
\begin{displaymath}
 \alpha_{1,p} \ldots  \alpha_{d,p}
\end{displaymath}
is well-defined as an element of $S^d (\fg^*)$.  Let
\begin{displaymath}
  j_p:pt \to M
\end{displaymath}
be the mapping ``$pt$'' onto $p$ and note that if $c$ is in
$H_G(M)$, $j^*_pc$ is in $H_G(pt)$ and thus in $S(\fg^*)$.  The
localization theorem asserts that, for every equivariant
cohomology class $c \in H_G(M)$,
\begin{equation}
  \label{eq:1.7}
  \int c = \sum_{p \in M^G}  \frac{j^*_pc}{\prod \alpha_{i,p}}
\end{equation}
%

There are many deep and beautiful applications of (\ref{eq:1.7})
but the focus of our interest in this article is that
(\ref{eq:1.7}) implies a lot of complicated identities among the
weights $\alpha_{i,p}$.  For instance, for $c=1$, it implies
\begin{equation}
  \label{eq:1.8}
  \sum (\prod \alpha_{i,p})^{-1} =0 \, .
\end{equation}
What are these identities?  
In particular, are there simpler identities of which they are 
formal consequences ?
We will show in part~2 of this article that one 
shed some light on these questions by looking at a
graph-theoretical analogue of
(\ref{eq:1.7}).

\subsection{The Jeffrey-Kirwan theorem}
\label{sec:1.3}

Another interesting source for identities of type (\ref{eq:1.7})
is the Jeffrey-Kirwan theorem (\cite{JK}): Let $K$ be a one-dimensional
connected closed subgroup of $G$ with Lie algebra $\fk$ and let
$\xi \in \fk$ be a basis vector of the group lattice of $K$. 
Suppose $M$ possesses a $G$-invariant symplectic form $\omega$ and 
that the action of $K$ on $M$ is Hamiltonian, i.e.
\begin{equation}
 \label{eq:1.9}
  \iota (\xi_M) \omega = - \; df \;\;,
\end{equation}
$f$ being a $G$-invariant function. In addition suppose that
$$M^K = M^G$$
and hence that the critical points of $f$ coincide with the
fixed points of $G$. Let $a$ be a regular value of $f$ and let
$M_a=f^{-1}(a)$. By the remark above, $M_a$ contains no $K$-fixed
points, so the action of $K$ on $M_a$ is locally free and the 
quotient space
$$M_a/K =  : M_{red}$$
is an orbifold. Moreover, from the action of $G$ on $M_a$ one 
gets an inherited action of the quotient group
$$ G/K =  : G_1$$
on $M_{red}$. Let $j$ be the inclusion of $M_a$ into $M$ and
$\pi$ the projection of $M_a$ onto $M_{red}$. By the 
Marsden-Weinstein theorem there exists a symplectic form 
$\omega_{red}$ on $M_{red}$ satisfying
$\pi^*\omega_{red}=j^*\omega$. In particular, $M_{red}$ is
{\em oriented}. If $\fg_1$ is the Lie algebra of $G_1$, its 
vector space dual, $\fg_1^*$, can be identified with the 
annihilator $\fg_{\xi}^*$ of $\xi$ in $\fg$; so there is an
integration operator
\begin{equation}
  \label{eq:1.10}
 \int  :  H_{G_1} ( M_{red}) \to S( \fg_{\xi}^* ) \; .
\end{equation}
Also, since the action of $K$ on $M_a$ is locally free, the map
$\pi$ induces an isomorphism
$$ \pi^*  :  H_{G_1}(M_{red}) \to H_G(M_a)$$
so one gets a map 
$$(\pi^*)^{-1} j^* =  : {\K}$$
of $H_G(M)$ into $H_{G_1}(M_{red})$. The Jeffrey-Kirwan theorem
asserts that for every equivariant cohomology class 
$c \in H_G(M)$,
\begin{equation}
  \label{eq:1.11}
  \int {\K} (c) = \Res_{\xi} \left( 
     \sum_{f(p) > a} \frac{j_p^*c}{\prod \alpha_{j,p}} \right) \; , 
     \quad p \in M^G \; ,
\end{equation}
$\Res_{\xi}$ being the residue of the rational function
in brackets with respect to the ``$\xi$-coordinate'' on $\fg^*$,
the other coordinates being held fixed. 
(This residue can be defined intrinsically to be an element
of $S(\fg_{\xi}^*)$. See \S\ref{sec:2.4})

\subsection{The Smith problem}
\label{sec:1.4}

The Smith conjecture asserts that if $M^G$ consists of two
points, $p$ and $q$, the isotropy representation of $G$ at $p$ is 
isomorphic (as a representation over $\RR$) to the isotropy
representation of $G$ at $q$.  The first complete proof of this
theorem (for $G$ an arbitrary compact Lie group) is due to
Atiyah, Bott and Milnor  (See \cite{AB2}, Theorem~3.83). If the
cardinality of $M^G$ is greater than two, the question of  
how the isotropy representations of $G$ at distinct 
fixed points are related to each other is still open and is known 
as the ``Smith problem''.  In this section we will describe a few 
of the more obvious relations:
\begin{enumerate}
\item Relations of type $J$.\\

Suppose that $M$ admits a $G$-invariant almost-complex
structure.  Then for every $p \in M^G$, the isotropy
representations of $G$ on $T_p$ is a complex representation, so
the weights of this representation,
\begin{equation}
  \label{eq:1.12}
  \alpha_{i,p}, \quad i=1, \ldots , d
\end{equation}
are unambiguously defined (not just defined up to sign).

For every closed subgroup $H$ of $G$ let $\fh$ be the Lie algebra 
of $H$ and
\begin{displaymath}
  \rho_H : \fg^* \to \fh^*
\end{displaymath}
the transpose of the inclusion map $\fh \to \fg$.  Let $E$ be a
connected component of $M^H$ and $p$ and $q$ elements of $E^G$.
We claim that the weights (\ref{eq:1.12}) can be ordered so that
\begin{equation}
\label{eq:1.13}
  \rho_H \alpha_{i,p} = \rho_H \alpha_{i,q} \, .
\end{equation}

\begin{proof}
  Let $x$ be an arbitrary point of $E$ and consider the isotropy 
  representation of $H$ on its normal space to $E$ at $x$.  This
  representation is a complex representation, so the weights of
  this representation are unambiguously defined and can't vary as 
  $x$ varies in $E$.  Thus, in particular, they have to be the
  same at $p$ and at $q$, implying (\ref{eq:1.13}).
\end{proof}

\item  Relations of type $\omega$.\\

Assuming that $M$ admits a $G$-invariant almost-complex
structure, $J$, is equivalent to assuming that $M$ admits a
$G$-invariant ``almost-symplectic'' structure, i.e.,~a two-form,
$\omega$, which is everywhere of maximal rank.  Suppose that 
$\omega$
is actually a symplectic form and the action $\tau$ is
\emph{Hamiltonian}, or, in other words, that there exists a
moment map
\begin{displaymath}
  \Phi : M \to \fg^* \, .
\end{displaymath}
From the convexity theorem ( \cite{A}, \cite{GS}) one gets

\begin{theorem*}
  Let $\Delta$ be the set of regular values of $\Phi$ in $\Phi
(M)$.  Then $\Delta$ is a disjoint union of open convex
polytopes.  Moreover, the vertices of these polytopes are the
images of the fixed points $p \in M^G$ and the edges going out
of these vertices are pointing in the directions of the vectors 
$\alpha_{i,p}$.

\end{theorem*}

The one-skeleton, $\Gamma$, of this configuration is called the
\emph{moment graph} of $M$.  (See \cite{Gu}.)  It exhibits a lot
of relations among the $\alpha_{i,p}$'s which are probably not
much simpler than the relations (\ref{eq:1.7}) but have the
virtue of being of a more geometric character.

\end{enumerate}

\subsection{The Goresky-Kottwitz-MacPherson graph}
\label{sec:1.5}

We will assume from now on that $M$ admits a $G$-invariant
almost-complex structure. Thus, for very $p \in M^G$, the
weights, $\alpha_{i,p} \in \fg^*$, are unambiguously defined.  In 
addition we will assume:  if $i \neq j$, $\alpha_{i,p}$ and
$\alpha_{j,p}$ are linearly independent.  This ``GKM hypothesis'' 
has the following consequence:  Let $\fh = \fh_i$ be the
annihilator of $\alpha_{i,p}$ in $\fg$ and let $H$ be the
$(n-1)$-dimensional subtorus of $G$ with Lie algebra, $\fh$.
Clearly, $p \in M^G \subseteq M^H$.

\begin{proposition*}
  Let $E$ be the connected component of $M^H$ containing $p$.
 Then
  \begin{displaymath}
    E \cong \CC P^1 \cong S^2
  \end{displaymath}
and the action of $G$ on $E$ is the standard action of the
circle, $S^1$, on $S^2$ by ``rotation about the $z$-axis''.  In
particular $E$ contains just two $G$-fixed points (one of which
is $p$).

\end{proposition*}

\begin{proof}
  The tangent space to $E$ at $p$ is the $2$-dimensional subspace 
  of $T_p$ on which $G$ acts with weight $\alpha_{i,p}$ so $E$
  itself is $2$-dimensional.  Since $E$ is compact and the action 
  of $G/H$ is non-trivial, $E$ is diffeomorphic to $S^2$, and
  this action is the standard action of $S^1$ on $S^2$ by the
  Korn-Lichtenstein theorem.

\end{proof}

Let $q$ be the other fixed point of $G$ in $E$.  Let
$\alpha_{p,e}= \alpha_{i,p}$ be the weight of the isotropy
representation of $G$ on $T_pE$ and let $\alpha_{q,e}=
\alpha_{j,q}$ be the weight of the isotropy representation of $G$
on $T_qE$.  From the fact that the action of $G$ on $E$ is
diffeomorphic to the standard action of $S^1$ on $S^2$ it follows 
that
\begin{equation}
  \label{eq:1.14}
  \alpha_{p,e}=- \alpha_{q,e} \, .
\end{equation}

For each of the weights $\alpha_{i,p}$ one gets an embedded
$\CC P^1$ of the type above, and we can represent these $\CC P^1$'s 
graphically by $d$ lines issuing from $p$.  Each of these lines
joins $p$ to another fixed point, $q$, and the $\CC P^1$'s
associated with the weights, $\alpha_{j,q}$, can also be
represented graphically by $d$ lines issuing from $q$.  One of
these will be the line from $p$ to $q$, but the remaining $d-1$
lines will join $q$ to other fixed points.  By repeating this
construction over and over until one runs out of fixed points,
one obtains a finite $d$-valent graph, $\Gamma$, the vertices of
which are the fixed points of $G$ and the edges of which
correspond to embedded $\CC P^1$'s, each of these $\CC P^1$'s
being a connected component of the fixed point set of an
$(n-1)$-dimensional subtorus of $G$.  We will call $\Gamma$ the
Goresky-Kottwitz-MacPherson (GKM)  graph of $M$.

\begin{example*}
  Suppose $M$ is a Hamiltonian $G$-manifold whose moment map,
  $\Phi$, maps $M^G$ injectively into $\fg^*$.  Then $\Phi$
  embeds the GKM graph into $\fg^*$, and its image is the moment
  graph.
\end{example*}

\subsection{The Goresky-Kottwitz-MacPherson theorem}
\label{sec:1.6}

The graph $\Gamma$ is equipped with an additional piece of
structure.  Namely let $I_{\Gamma}$ be the incidence relation of
this graph:  the set of all pairs $(p,e)$, $p$ being a vertex and
$e$ an edge containing $p$.  Then one has a  map
\begin{equation}
  \label{eq:1.15}
  \alpha : I_{\Gamma} \to \fg^* - \{ 0 \}
\end{equation}
mapping $(p,e)$ to the weight $\alpha_{p,e}$.  We will call this 
the \emph{axial function} of $\Gamma$.  It has the following
properties (the first two of which we have already commented on):

\begin{enumerate}
\item If $e$ is an edge and $p$ and $q$ are the vertices joined
  by $e$
  \begin{equation}
    \label{eq:1.16}
      \alpha_{p,e} = - \alpha_{q,e}
  \end{equation}

\item  If $p$ is a vertex and $e_1, \ldots , e_d$ are the edges
  containing $p$, the vectors
  \begin{equation}
    \label{eq:1.17}
    \alpha_{p,e_i} , \quad i=1, \ldots ,d
  \end{equation}

are pair-wise linearly independent.

\item  Let $e$ be an edge and $p$ and $q$ the vertices joined by
  $e$.  Let
  \begin{displaymath}
    \fg_e = \{ \xi \in \fg , \alpha_{p,e}(\xi)=0 \}
  \end{displaymath}
and let
\begin{displaymath}
    \rho_e : \fg^* \to \fg^*_e
  \end{displaymath}
be the transpose of the inclusion map, $\fg_e \to \fg$.  Let
$e_i, i=1, \ldots,d$ and $e'_i , i=1, \ldots ,d$ be the edges
containing $p$ and $q$ respectively, with $e_d=e'_d=e$.  Then the
$e_i$'s can be ordered so that
\begin{equation}
  \label{eq:1.18}
  \rho_e \alpha_{p,e_i} = \rho_e \alpha_{q,e'_i} \, .
\end{equation}

\end{enumerate}

\begin{proof}
  (\ref{eq:1.18}) is just a special case of (\ref{eq:1.13}), $H$
  being the subtorus of $G$ with Lie algebra $\fg_e$.
\end{proof}

From the data $(\Gamma, \alpha)$ one can construct a graded ring
\begin{equation}
  \label{eq:1.19}
  H(\Gamma, \alpha)= \bigoplus H^{2k}(\Gamma, \alpha)
\end{equation}
as follows.  For each edge $e$, the map $\rho_e : \fg^* \to
\fg^*_e$ extends to a ring morphism 
\begin{equation}
  \label{eq:1.20}
  \rho_e : S (\fg^*) \to S (\fg^*_e) \, .
\end{equation}
Let $V_{\Gamma}$ be the set of vertices of $\Gamma$ and let
$H^{2k}(\Gamma, \alpha)$ be the set of all maps
\begin{displaymath}
  f: V_{\Gamma} \to S^k (\fg^*)
\end{displaymath}
satisfying the compatibility condition:
\begin{equation}
  \label{eq:1.21}
  \rho_e f(p) = \rho_e f (q)
\end{equation}
for all vertices, $p$ and $q$, and edges, $e$, joining $p$ to
$q$.  $H(\Gamma, \alpha)$ can be given a ring structure by
pointwise multiplication
\begin{displaymath}
  (f_1f_2) (p) = f_1 (p) f_2(p) \, .
\end{displaymath}
(Notice that if $f_1$ and $f_2$ satisfy (\ref{eq:1.21}) so does
$f_1f_2$ since $\rho_e$ is a ring morphism.)  In addition,
$H(\Gamma, \alpha)$ contains $S(\fg^*)$ as a subring:  the ring
of constant maps of $V_{\Gamma}$ into $S(\fg^*)$.  In particular, 
$H(\Gamma, \alpha)$ is a module over $S(\fg^*)$.

\begin{theorem*}(\cite{GKM})
If $M$ is equivariantly formal, $H_G (M)$ is isomorphic, as a
graded ring, to $H(\Gamma, \alpha)$.  
\end{theorem*}


\subsection{Holonomy}
\label{sec:1.7}

Let $\PP$ be the complex projective line, let $z$ be the standard 
coordinate function on $\CC = \PP_0 = \PP- \{ \infty \}$ and
$w=z^{-1}$ the corresponding coordinate function on
$\PP_{\infty}= \PP - \{ 0 \}$.  The multiplicative group,
$\CC^{\#}= \CC - \{ 0 \}$, acts on $\CC$ by homotheties:
\begin{displaymath}
  (a,z) \to az
\end{displaymath}
and this extends to a holomorphic action, $\rho$, of $\CC^\#$ on
$\PP$.  Let $E$ be a holomorphic, rank $r$ vector bundle over
$\PP$ and suppose $\rho$ lifts to an action, $\rho$, of $\CC^\#$
on $E$ by vector bundle automorphisms.  Let $E_0$ and
$E_{\infty}$ be fibers of $E$ over $0$ and $\infty$.  The
isotropy representations of $\CC^\#$ at $0$ and $\infty$
decompose these spaces into invariant one-dimensional subspaces
\begin{equation}
  \label{eq:1.22}
  E_0 = V_1 \oplus \cdots \oplus V_r
\end{equation}
and 
\begin{equation}
  \label{eq:1.23}
  E_{\infty} = V'_1 \oplus \cdots \oplus V'_r \, .
\end{equation}
Let $m_1 , \ldots , m_r$ be the weights of the representations of 
$\CC^\#$ on $V_1 , \ldots , V_r$.  We will assume that these
weights are all distinct and that $m_1 < m_2 < \cdots < m_r$.
Similarly we will assume that the weights $m'_i$, of the
representations of $\CC^\#$ on $V'_1 , \ldots , V'_r$ are all
distinct (but we won't require that $m'_1 < \cdots < m'_r$).  By
a theorem of Birkhoff-Grothendieck (\cite{Ha}) 
there is an equivariant decomposition of $E$ into line-bundles
\begin{equation}
  \label{eq:1.24}
  E=\LL_1 \oplus \cdots \oplus \LL_r
\end{equation}
such that the fiber of $\LL_i$ over $0$ is $V_i$.  Moreover this
decomposition is \emph{unique up to isomorphism}.  (To see this,
let
\begin{displaymath}
  E= \LL'_1 \oplus \cdots \oplus \LL'_r
\end{displaymath}
be another decomposition with these properties.  Over $\PP_0$,
one can find a trivializing section, $s_i$, of $\LL_i$ which
transforms under $\CC^\#$ according to the law
\begin{equation}
  \label{eq:1.25}
  \rho^*_a s_i = a^{m_i} s_i
\end{equation}
and a trivializing section, $s'_i$, of $\LL'_i$ with the same
transformation property.  Moreover, one can assume that $s_i(0) = 
s'_i (0)$.  We claim that $s_1 \equiv s'_1$ (and hence $\LL_1 = 
\LL'_1$).  To see this we note that because of the transformation 
law (\ref{eq:1.25})
\begin{displaymath}
  s'_1 = s_1 + \sum_{i>1} c_i z^{m_1-m_i} s_i
\end{displaymath}
on $\PP_0$.   However, since $s'_1$ is holomorphic near $0$ and
$m_1 - m_i$ is strictly negative, the constants, $c_i$, are all
zero.

Applying this argument to the quotient bundle, $E/\LL_1$, one
concludes by induction that $\LL_i \cong \LL'_i$ for all $i$. Q.E.D.)

In particular, coming back to the isotropy decompositions
(\ref{eq:1.22})  and (\ref{eq:1.23}), one has a canonical map of
the $r$-element set
\begin{equation}
  \label{eq:1.26}
  \{ V_1 , \ldots , V_r \}
\end{equation}
onto the $r$-element set
\begin{equation}
  \label{eq:1.27}
  \{ V'_1 , \ldots , V'_r \}
\end{equation}
which maps $V_i = (\LL_i)_0$ onto $(\LL_i)_{\infty}$; and we can
relabel the $V'_i$'s so that this map maps the element, $V_i$, of 
the set (\ref{eq:1.26}) onto the element, $V'_i$, of the set
(\ref{eq:1.27}).  We will call the map the \emph{holonomy} map.

Let us apply these remarks to one of the $\CC P^1$'s in
\S\ref{sec:1.5}; i.e.,~let $H$ be a subtorus of $G$ of
codimension one and $X=\CC P^1$ a $2$-dimensional connected
component of $M^H$.  Let $NX$ be the normal bundle of $X$ and let 
$\beta_1 , \ldots , \beta_N$ be the weights of the isotropy
representation of $H$ on a typical fiber of $X$.  The fiber of
$N$ over $x \in X$ splits into a direct sum of vector subspaces
\begin{equation}
\label{eq:1.28}
  (E_1)_x \oplus \cdots \oplus (E_N)_x
\end{equation}
such that the weight of the representation of $H$ on $(E_i)_x$ is 
$\beta_i$.  Since the $\beta_i$'s can occur with multiplicities,
the dimension of  $(E_i)_x$ doesn't have to be one; however the
assignment, $x \to (E_i)_x$, defines a complex vector bundle,
$E_i$, over $X$.  Moreover, this bundle can be given a
holomorphic structure, and the action of $S^1=G/H$ on $E_i$ can
be extended to a holomorphic action of $\CC^\#$.  In particular
one can define for $E_i$ a \emph{holonomy} map of the type we
described at the beginning of this section.

Now let $X^G= \{ p,q \}$, and $\alpha_{i,p}$ and $\alpha_{i,q}$
be the weights of the isotropy representations of $G$ on $T_pM$
and $T_qM$.  By assumption these weights occur with multiplicity
one, so one gets decompositions of $T_pM$ and $T_qM$ into direct
sums of one-dimensional weight spaces
\begin{equation*}
  T^{\alpha_{1,p}}_p \oplus \cdots \oplus T^{\alpha_{d,p}}_p 
\end{equation*}
and
\begin{equation*}
  T^{\alpha_{1,q}}_q \oplus \cdots \oplus T^{\alpha_{d,q}}_q \, .
\end{equation*}
Let $\rho $ be the projection of $\fg^*$ onto $\fh^*$.  By
(\ref{eq:1.18}) we can reorder the $\alpha_{i,q}$'s so that $\rho 
(\alpha_{i,p})= \beta_i= \rho (\alpha_{i,q})$; so, if the
$\beta_i$'s are all distinct, one has a \emph{canonical} map of
the set
\begin{equation}
  \label{eq:1.29}
  \{ \alpha_{1,p}, \ldots , \alpha_{d,p} \}
\end{equation}
onto the set
\begin{equation}
  \label{eq:1.30}
  \{ \alpha_{1,q}, \ldots , \alpha_{d,q} \}
\end{equation}
mapping $\alpha_{i,p}$ onto $\alpha_{i,q}$.  However, this
canonical map can even be defined when the $\beta_i$'s are
\emph{not} distinct (that is, when the vector bundles, $E_i$, are 
\emph{not} of rank one) by using the holonomy structure on these
bundles.  In other words there exists a canonical holonomy
mapping from the set (\ref{eq:1.29}) to the set (\ref{eq:1.30})
which, when the $\beta_i$'s are distinct, is defined trivially by 
the recipe
\begin{displaymath}
  \alpha_{i,p} \to \rho (\alpha_{i,p}) = \beta_i 
  = \rho (\alpha_{i,q}) \to \alpha_{i,q}
\end{displaymath}
but, when the $\beta_i$'s are not distinct, involves some
topological properties of the bundles, $E_i$.  We will denote
this map by $\theta_{p,e}$.

One can give a slightly more ``graphical'' description of this
map:  As in \S\ref{sec:1.6}, let $\Gamma$ be the GKM graph, let $
V_{\Gamma}$ be its vertices and let $I_{\Gamma}$ be the incidence 
relation.  The projection, $I_{\Gamma} \to V_{\Gamma}$, can be
regarded as a fiber bundle over $V_{\Gamma}$, the fiber, $E_p$,
over $p \in V_{\Gamma}$ being the set of all pairs, $(p,e)$ in
$I_{\Gamma}$ (i.e.,~the set of all oriented edges of $\Gamma$
pointing out from $p$).  In this ``fiber bundle'' picture,
$\theta_{p,e}$ is just a map
\begin{equation}
  \label{eq:1.31}
  \theta_{p,e} : E_p \to E_q
\end{equation}
with the properties 
\begin{equation}
  \label{eq:1.32}
  \theta_{p,e} (p,e) = (q,e)
\end{equation}
and
\begin{equation}
  \label{eq:1.33}
  \theta_{p,e} = \theta^{-1}_{q,e}
\end{equation}
$e$ being the edge joining $p$ to $q$.  A family of maps,
$\theta$, with the properties (\ref{eq:1.32})--(\ref{eq:1.33}) is 
called a \emph{connection} on $\Gamma$ (cf.~\cite{CS}).  In terms 
of this connection, one can reformulate (\ref{eq:1.18}) more
precisely:

\begin{theorem}
  
  Let $e$ be an edge of $\Gamma$ joining the vertex, $p$, to the
  vertex, $q$.  Then for every $(p,e_i) \in E_p$
  \begin{equation}
    \label{eq:1.34}
    \rho_e (\alpha_{p,e_i}) = \rho_e (\alpha_{q,e'_i})
  \end{equation}
$(q,e'_i)$ being the image with respect to $\theta_{p,e}$ of $(p,e_i)$.

\end{theorem}

\subsection{Betti numbers}
\label{sec:1.8}

The theorem of Goresky-Kottwitz-MacPherson described in \S~\ref{sec:1.6} implies that the odd Betti numbers of $M$ are 
zero. The even Betti numbers can be computed as follows.
As in \S\ref{sec:1.3} let $K$ be a one-dimensional closed connected 
subgroup of $G$ with $M^G=M^K$ and let $\xi$ be a basis vector of
$\fk$. For every $p \in V_{\Gamma}$ let $\sigma_p$ be the number of 
edges, $e$, with one vertex $p$ and with $\alpha_{p,e}(\xi) <0$.

\begin{theorem*}
The $2k$-th Betti number, $\beta_k$, of $M$ is equal to the number of
points, $p \in V_{\Gamma}$, with $\sigma_p = k$.
\end{theorem*}

If $M$ possesses a $G$-invariant symplectic form having the 
properties described in \S\ref{sec:1.3}, this theorem can be 
proved by Morse theory : Let $f$ be the function defined by 
(\ref{eq:1.9}). The critical points of this function coincide
with the fixed points of $G$ and it is not difficult to show that 
the index of the Hessian of $f$ at $p \in M^G$ is just $2\sigma_p$.

Note by the way that since $M$ is equivariantly formal, the identity
(\ref{eq:1.6}), coupled with the Goresky-Kottwitz-MacPherson theorem,
implies that
\begin{equation}
\label{eq:1.35}
\mbox{dim }H^{2k}(\Gamma,\alpha) = \sum \beta_r 
\mbox{ dim }S^{k-r}(\fg^*).
\end{equation}
Thus the {\it additive} structure of $H(\Gamma,\alpha)$ can be computed 
just by inspecting the orientations of the edges of $\Gamma$.

\setcounter{section}{2}
\setcounter{subsection}{0}
\setcounter{equation}{0}
\subsection{Graphs and axial functions}
\label{sec:2.1}

In part~2 of this article, $\fg$ will simply be a vector space
over $\RR$ of dimension $n$ and $\fg^*$ its vector space dual.
(In particular $\fg$ will \emph{not} necessarily be the Lie
algebra of a group, $G$.)  Let $\Gamma$ be a finite simple $d$-valent graph and
$I_{\Gamma}$ its incidence relation.

\begin{definition}{1.}
\label{def:1}
An {\em axial function} on $\Gamma$ is a map $\alpha: I_{\Gamma}
\to \fg^* - \{ 0 \}$ having the properties
(\ref{eq:1.16})--(\ref{eq:1.18}).
\end{definition}

Note that by axiom (\ref{eq:1.16}) an axial function can be
thought as a function
$$\alpha : E_{\Gamma}^{\pm} \to \fg^* $$
on the set of {\it oriented} edges of $\Gamma$, having the 
property that if $e^+$ and $e^-$ are the two oriented edges of 
$\Gamma$ associated with an unoriented edge $e$, then 
$$\alpha (e^+) =  -\alpha (e^-).$$
\begin{definition}
\label{def:2}
A \emph{connection}, $\theta$, on $\Gamma$ is a collection of maps,
$\theta_{p,e}$, $(p,e) \in I_{\Gamma}$, satisfying the axioms
(\ref{eq:1.31})--(\ref{eq:1.33}); an axial function,
$\alpha$, is \emph{compatible} with $\theta$ if it satisfies the
stronger version, (\ref{eq:1.34}) of axiom~(\ref{eq:1.18}).

\end{definition}

For instance, suppose that for every vertex, $p$, and every edge, 
$e$, containing $p$, no two of the vectors (\ref{eq:1.18}) are
equal.  Then, as we showed in section~\ref{sec:1.7}, there is a
  \emph{unique} connection on $\Gamma$ which is compatible with $ 
 \alpha$.

We will list below a few examples of axial functions and
connections on graphs and describe some of their functorial properties:

\begin{example}{\it The complete graph on $N$ vertices. }
The vertices of this graph are the elements of the $N$-element set
$\{ 1,\dots , N \}$ and each pair of elements, $(i,j),\; i \neq j$,
is joined by an edge. Thus the set of oriented edges is just the
set
$$ \{ (i,j) ; 1 \leq i,j \leq N, \; i \neq j \}.$$
Let $\alpha_i, \; i=1,\dots, N$ be non-zero elements of $\fg^*$.
Then the function $(i,j) \to \alpha_i - \alpha_j$ satisfies 
(\ref{eq:1.16}) and (\ref{eq:1.18}) and hence it is  an axial 
function iff, for every $i$, the $N-1$ vectors $\alpha_i -\alpha_j, 
\; i\neq j$ are pairwise linearly independent. Conversely, we claim that
{\it every} axial function is of this form. 
(Proof : Let $(i,j) \to \alpha_{i,j}$ be an axial function. Then if 
$i,j$ and $k$ are distinct
$$\alpha_{i,k}=\alpha_{j,k} + c_{i,j} \alpha_{i,j} 
\qquad c_{i,j}=c_{j,i}.$$
Hence
$$\alpha_{k,i}=\alpha_{j,i} + c_{k,j} \alpha_{k,j}=
\alpha_{k,j}+c_{i,j} \alpha_{j,i}$$
so $(1-c_{i,j})\alpha_{j,i} = ( 1-c_{k,j})\alpha_{k,j}$ and
$c_{i,j}=c_{k,j}=1$. Now let $\alpha_1=0$ and let 
$\alpha_i =\alpha_{i,1}$ for $i>1$.)
\end{example}

An $\alpha$-compatible connection, $\theta$, can be defined as
follows:  For every oriented edge, $(i,j)$, of $\Gamma$ let
$\theta_{i,j}$ be the map of the set of edges of $\Gamma$
containing $i$ onto the set of edges of $\Gamma$ containing $j$
which maps $(i,j)$ onto $(j,i)$ and, for $k \ne i,j$, maps
$(i,k)$ onto $(j,k)$.

\begin{example}{\it Sub-objects. } \label{ex:2}
Let $\Gamma_1$ be an $r$-valent sub-graph of $\Gamma$ and $j$ the
embedding of $V_{\Gamma_1}$ into $V_{\Gamma}$. From $j$ one gets 
an embedding  $i : I_{\Gamma_1} \to I_{\Gamma}$ and one can 
pull-back the axial function $\alpha$ to $I_{\Gamma_1}$. In 
general $i^*\alpha$ won't be an axial function; however if it is,
we will say that $\Gamma_1$ is {\it compatible} with $\alpha$ and
call $(\Gamma_1, i^*\alpha)$ a {\it sub-object} of $(\Gamma, \alpha)$.
\end{example}


\begin{example}
  \label{ex:3}
\emph{Sub-objects which are ``totally geodesic'' with respect to a connection.}
  Let $\theta$ be a connection on $\Gamma$, and let $\Gamma'$ be
  an $r$-valent subgraph of $\Gamma$.  For every vertex, $p$, of
  $\Gamma$ let $E_p$ and $E'_p$ be the oriented edges of $\Gamma$ 
  and of $\Gamma'$ issuing from $p$.  We will say that $\Gamma'$
  is totally geodesic with respect to $\theta$ if, for every
  oriented edge, $(p,q)$, of $\Gamma'$ the restriction of the
  holonomy map, $\theta_{p,q} \, : \, E_p \to E_q$, to $E'_p$
  maps $E'_p$ onto $E'_q$.

If this happens, this restriction defines an induced connection,
$\theta'$, on $\Gamma'$.  Moreover, if $\alpha : I_{\Gamma} \to
\fg^*$ is an axial function which is compatible with $\theta$,
the restriction of $\alpha$ to $I_{\Gamma'}$ is an axial function 
and is compatible with $\theta'$.
\end{example}

\begin{example}
    \label{ex:4}
\emph{The totally geodesic sub-objects of $\Gamma_N$.}
Let $\Gamma_N$ be the complete graph on $N$ vertices and for
every subset, $S$, of $\{1, \ldots ,N \}$ let $\Gamma_S$ be the
graph whose  vertices are the elements of $S$ and whose oriented
edges are the pairs $(s_1,s_2)$, $(s_i \in S), s_1 \ne s_2$.  It
is obvious that $\Gamma_S$ is a totally geodesic with respect to
the connection we defined in example~1; and, in fact,
every totally geodesic sub-object of  $\Gamma_N$ is a $\Gamma_S$
for some subset $S$. 
(To see this, let $\Gamma'$ be a
connected totally geodesic subgraph and let $(p_1,p_2)$ be an
oriented edge of $\Gamma'$.  
If $q$ is a vertex of $\Gamma'$ distinct from $p_1$ and $p_2$ and 
$(p_1,q)$ is an oriented edge of $\Gamma'$, $(p_2,q)$ has to be an
oriented edge of $\Gamma'$,so it follows from the
connectivity of $\Gamma'$ that, for every pair of vertices, $p$
and $q$, of $\Gamma'$ $(p,q) $ is an oriented edge of $\Gamma'$.)
\end{example}

\begin{example}{\it The graph $\Gamma_{\fh}$.}
\label{ex:5}
Let $\fh$ a vector
subspace of $\fg$ and let $\rho_{\fh}~:~ \fg^*~ \to~ \fh^*$ be 
the transpose of the inclusion map $\fh \to \fg$. 
Let $\Gamma_{\fh}$ be the subgraph of $\Gamma$ whose edges are 
the edges, $e$, of $\Gamma$ for which 
\begin{equation}
\label{eq:2.1}
\rho_{\fh} \alpha_{p,e} = - \rho_{\fh} \alpha_{q,e} =0
\end{equation}
$p$ and $q$ being the vertices of $e$. Each connected component of
this graph is $k$-valent for some $k$ and is a sub-object in the sense
of item~2. Moreover, if $p$ and $q$ are in the same 
connected component of $\Gamma_{\fh}$ and $e_i$ and 
$e'_i$, $i=1,...,d$, are the edges of $\Gamma$ containing $p$ and $q$,
one can order the $e_i$'s so that 
\begin{equation}
\label{eq:2.2}
\rho_{\fh} \alpha_{p,e_i} =\rho_{\fh} \alpha_{q,e'_i}
\end{equation}
(compare with (\ref{eq:1.13})).
\end{example}

\begin{example}{\it Product objects. } 
\label{ex:6}
Let $\Gamma_1$ be a graph of 
valence $d_1$ and $\Gamma_2$ a graph of valence $d_2$. The vertices
of the {\it product} graph, $\Gamma_1 \times \Gamma_2$, are the 
pairs $(p,q)$, $p \in V_{\Gamma_1}$ and $q \in V_{\Gamma_2}$; 
two vertices $(p,q)$ and $(p',q')$ are joined by an edge if either
$p=p'$ and $q$ and $q'$ are joined by an edge in $\Gamma_2$ or
$q=q'$ and $p$ and $p'$ are joined by an edge in $\Gamma_1$. Thus 
this product graph is a $d_1+d_2$-valent graph and its set of oriented 
edges is the disjoint union
\begin{equation}
\label{eq:2.3}
E_{\Gamma_1}^{\pm} \amalg E_{\Gamma_2}^{\pm}
\end{equation}
of the set of oriented edges of $\Gamma_1$ and the set of oriented edges
of $\Gamma_2$. If
$$\alpha_i : E_{\Gamma_i}^{\pm} \to \fg^*, \qquad i=1,2$$
is an axial function on $\Gamma_i$, one defines the 
product function $\alpha$ on $\Gamma_1 \times \Gamma_2$ to be
the function which is equal to $\alpha_1$ on the first summand
of (\ref{eq:2.3}) and equal to $\alpha_2$ on the second
summand. 
Then $\alpha$ satisfies (\ref{eq:1.16}) - (\ref{eq:1.18}) and 
it is called 
{\em the product axial function}.
If in addition, $\theta_1$ and $\theta_2$ are
connections on $\Gamma_1$ and $\Gamma_2$ one can define a product
connection, $\theta$, on $\Gamma_1 \times \Gamma_2$ by letting
\begin{displaymath}
  \theta_{(p,q;p',q)}= \theta_{p,p'} \times (\Identity)_{q,q}
\end{displaymath}
and letting
\begin{displaymath}
  \theta_{(p,q;p,q')} = (\Identity)_{p,p} \times \theta_{q,q'}
  \, .
\end{displaymath}
If $\alpha_i$ and $\theta_i$ are compatible, for $i=1$ and $2$,
the product axial function which we defined above is compatible
with $\theta$.

\end{example}

\begin{example}{\it Blowing-up. } 
\label{ex:7}
This operation can be defined for 
any sub-object of a graph-axial function pair; however, 
for simplicity, we will only consider here the special case when 
the sub-object is a point. Let $\Gamma$ be a finite simple $d$-valent 
graph and let $p_0$ be an arbitrary vertex of $\Gamma$. Let 
$e_i,\; 1=1,...,d$ be the edges of $\Gamma$ containing $p_0$ and let 
$q_i$ be the vertex joined by $e_i$ to $p_0$. From this data one can
 construct a new graph, $\Gamma^{\#}$, as follows. Replace the vertex 
$p_0$ by $d$ new vertices, $p_i,\; i=1,...,d$ ( which one should 
think of as being the ``baricenters'' of the edges $e_i$) and to
 each of these new vertices adjoin $d$ edges; one edge going from
$p_i$ to $q_i$ (which one should think of as being a replacement for 
the old edge $e_i$) and one edge going from $p_i$ to each of the $p_j's$,
$j \neq i$. Let 
\begin{equation}
\label{eq:2.4}
\beta : V_{\Gamma^{\#}} \to V_{\Gamma}
\end{equation}
be the map which sends $\{ p_1,...,p_d \}$ to $p_0$ and is  the 
identity on the complement of $\{ p_1,...,p_d \}$. We will call this 
map the {\it blowing-down} map. The set $\{ p_1,...,p_d\}$, which is
the pre-image of $p_0$ with respect to $\beta$, is the set of vertices 
of a sub-graph, $\Gamma_0$, of $\Gamma^{\#}$ (the complete graph on
$d$ vertices) which we will call {\it the singular locus} of the
blowing-down map $\beta$.

Now let $\alpha : I_{\Gamma} \to \fg^*$ be an axial function and let
$\alpha_i, \; i=1,...,d$ be the values of $\alpha$ on the vertex-edge 
pairs $(p_0, e_i)$. Let us assume that for each $i$ the $d-1$ vectors 
$\alpha_j - \alpha_i, \; j \neq i$ are pairwise linearly independent. 
We can then define an axial function, $\alpha^{\#}$, on $\Gamma^{\#}$,
as follows:
\begin{enumerate}
\item On the oriented edges, $e$, of $\Gamma$, not containing $p_0$,
        $\alpha^{\#} (e)= \alpha (e)$.
\item On the oriented edges, $e=(p_i, q_i)$, 
        $\alpha^{\#} (e)= \alpha_i$.
\item On the oriented edges, $e=(p_i,p_j)$, 
        $\alpha^{\#} (e)= \alpha_j-\alpha_i$.
\end{enumerate}

This defines $\alpha^{\#}$ on {\it all} edges of $\Gamma^{\#}$ and it is
easy to check that $\alpha^{\#}$ satisfies the axioms 
(\ref{eq:1.16}) - (\ref{eq:1.18}).

Finally, if $\theta$ is a connection on $\Gamma$ there is a
unique connection, $\theta^\#$, on $\Gamma^\#$ with the following
properties.

\begin{enumerate}
\item The restriction of $\theta$ to $\Gamma_0$ is the connection 
  described in example~one.

\item  $\theta_{(p_i,p_j)}^\#$ maps the oriented edge,
  $(p_i,q_i)$, to the oriented edge $(p_j,q_j)$.

\item  $\theta^\#_{(p_i,q_i)}$ maps the oriented edge $(p_i,q_i)$ 
  to the oriented edge $\theta_{(p,q_i)} (e_j)$.

\item  If $q$ is not equal to one of the $p_i$'s
 and $q$ and $q_i$ are joined by an edge, 
\begin{displaymath}
\theta^\#_{(q_i,q)}
 (e^\#_i) = \theta_{(q_i,q)} (e^-_i) \, ,
\end{displaymath}
 $e^\#_i$ being the edge
 joining $q_i$ to $p_i$ and $e^-_i$ the edge joining $q_i$ to
 $p$.  On the other edges of $\Gamma^\#$ issuing from $q_i$,
 $\theta^\#_{(q_i,q)}= \theta_{(q_i,q)}$.

\item  If $q$ and $q'$ are not equal to one of the $q_i$'s or one 
  of the $p_i$'s, $\theta^\#_{(q,q')} = \theta_{(q,q')}  $.

\end{enumerate}

\end{example}

\begin{example}{\it  The case $d$=dim  $\fg^*$ =2 . }
\label{ex:8}
Let $\Gamma$ be a finite connected 2-valent graph with $N$ vertices, 
$\fg^*$ a 2-dimensional vector space and $\alpha$ an axial function. 
We will orient the edges of $\Gamma$ so that for each vertex, $p$, one
of the edges containing $p$ is pointing in the direction of $p$ and the 
other is pointing away from $p$ (since $\Gamma$ is connected there 
are clearly just two ways of orienting the edges so that their 
orientations have this property.). Let $p_1,...,p_{N},p_{N+1}=p_1$ be
an enumeration of the vertices of $\Gamma$ such that the outward 
pointing edge at $p_i$ joins $p_i$ to $p_{i+1}$ and let $\alpha_i$ be
the value of $\alpha$ on $(p_i,p_{i+1})$. Then (\ref{eq:1.18}) is 
equivalent to
\begin{equation}
\label{eq:2.5}
\alpha_{i+1} \wedge \alpha_i = \alpha_i \wedge \alpha_{i-1}
\end{equation}
for all $i$. ( For example for $N=4k$, let $ \{ \alpha_1, \alpha_2 \}$
be a basis of $\fg^*$. Then a solution of (\ref{eq:2.5}) is obtained
by letting $\alpha_1= - \alpha_3 = \alpha_5 = ...$ and 
$ \alpha_2 = - \alpha_4 = \alpha_6 = ...$).
\end{example}

\subsection{Orientations}
\label{sec:2.2new}

Let $(\Gamma, \alpha)$ be a graph-axial function pair and let
\begin{displaymath}
\P = \{ \xi \in \fg , \alpha_{p,e} (\xi) \neq 0
\hbox{ for all } (p,e) \in I_{\Gamma} \} \, .
\end{displaymath}
Then for every $\xi \in \P$, the axial function $\alpha$
defines an \emph{orientation} of $\Gamma$; in other words, for
each edge $e$, it fixes an ordering of the vertices of $e$.
Namely if $p$ and $q$ are the vertices of $e$, one orders them so 
that
\begin{equation}
  \label{eq:2.6}
p< q \Leftrightarrow \alpha_{p,e} (\xi) >0 \, .
\end{equation}
It is clear that this orientation doesn't depend on $\xi$ but only
on the connected component of $\P$ in which $\xi$ is contained.
On the other hand it is clear that different components will give
rise to different orientations (for instance, replacing $\xi$ by
$-\xi$ reverses all the orientations). We will say that $\Gamma$
satisfies the {\em no-cycle condition} if, for at least one of these
orientations, $\Gamma$ has no cycles.

\begin{definition} 
Given $\xi \in \P$, a function $f: V_{\Gamma} \to \RR$ is
\emph{positively oriented with respect to} $\xi$ if, for every
pair of vertices, $p$ and $q$, and edge, $e$, joining $p$ to $q$,
the ratio of $f(p) - f(q)$ to $\alpha_{q,e}(\xi)$ is positive.

\end{definition}

If $f$ is positively oriented with respect to $\xi$, the
orientation of $\Gamma$ associated with $\xi$ can't have closed
 cycles since $f$ has to be strictly increasing along any
oriented path.  We will prove that the converse is true:

\begin{theorem*}
  If the orientation of $\Gamma $ associated with $\xi$ has no
  cycles, there exists a function $f: V_{\Gamma} \to \RR$ which 
  is positively oriented with respect to $\xi$.
\end{theorem*}

\begin{proof}
  Given $p \in V_{\Gamma}$, consider the \emph{longest} oriented
  path with initial point $p$, i.e.,~the longest sequence
  \begin{equation}
    \label{eq:2.7}
    p=p_1,p_2, \ldots , \quad p_i \in V_{\Gamma}
  \end{equation}
with the property that $p_i$ and $p_{i+1}$ are the vertices of a
common edge and, relative to the orientation on this edge, $p_i < 
p_{i+1}$.  If $\Gamma$ has no cycles this longest path has to be of
finite length, i.e.,~has to terminate at some point, $p_N$.  Now
set $f(p)=-N$.  It is easy to check that this function is
positively oriented with respect to $\xi$.

\end{proof}

{\bf Remarks:}
  \begin{enumerate}
  \item The vertices, $p$, where $f(p)=0$ have the property that
    all edges containing $p$ are pointing ``into'' $p$, i.e.~$p$ 
    is a ``maximum'' of the oriented graph $\Gamma$.  In
    particular, if $f(p)=-N$, this is true of the vertex $p_N$
    in the sequence (\ref{eq:2.7}); so the argument above shows
    that every vertex can be joined by an oriented path to a
    maximal vertex.

\item  One can perturb $f$ so that it remains positively oriented 
  with respect to $\xi$ and, in addition, takes on distinct
  values at distinct vertices.  Namely suppose that 
  \begin{displaymath}
    f^{-1} (k) = \{ p_1 , \ldots , p_r \} \, .
  \end{displaymath}
Redefine $f$ on the set $\{ p_1 , \ldots , p_r \}$ by setting $f$ 
equal to $k+ \epsilon_i$ on $p_i$ where $\epsilon_i \neq
\epsilon_j$ for $i \neq j$ and the $\epsilon_i$'s are small.
This redefined function is still positively oriented with respect 
to $\xi$ but now takes distinct values at $p_1 , \ldots , p_r$.
\end{enumerate}

\subsection{The cohomology ring of $(\Gamma, \alpha)$}
\label{sec:2.2}

We define the \emph{cohomology ring} of $(\Gamma, \alpha)$ to be 
the ring $H(\Gamma, \alpha)$ which we defined in 
\S\ref{sec:1.6}.  
As $(\Gamma, \alpha)$ is no longer the GKM
data associated with a $G$~-~manifold, it is, perhaps, a misnomer
to refer to this ring as a ``cohomology ring''; however, there
are other reasons for using this terminology.  For instance, if $
\Gamma$ is the one-skeleton of a simplicial polytope, 
$H(\Gamma,  \alpha)$ is just the 
Stanley-Reissner cohomology ring of the
dual polytope.  (We are indebted to Mark Goresky for this
observation.)  We will describe below a few properties of this
ring.

\begin{enumerate}
\item \label{item:1} As we pointed out in \S\ref{sec:1.6}, 
 $H(\Gamma,\alpha)$ contains $S (\fg^*)$ as a subring.

\item  \emph{Chern classes:}  For each $p \in V_{\Gamma}$, let
  $e_1, \ldots ,e_d$ be the edges containing $p$ and let $c_k
  (p)$ be the $k$-th elementary symmetric function  in the
  monomials $\alpha_{p,e_1}, \ldots , \alpha_{p, e_d}$.  The
  function $p \to c_k (p)$ defines an element $c_k$ of
  $H^{2k}(\Gamma, \alpha)$ which can be thought of as the $k$-th 
  Chern class of the ``tangent bundle'' of $\Gamma$.

\item  \emph{Symplectic structures:}  An element of 
$H^2(\Gamma,\alpha)$ is just a map 
$c: V_{\Gamma} \to \fg^*$ satisfying
  \begin{equation}
    \label{eq:2.8}
    c(p) - c(q) = \lambda_e \alpha_{q,e}
  \end{equation}
for every pair of vertices, $p,q$, and edge, $e$, joining $p$
to $q$.  We will call $c$ {\em symplectic} if, for every edge $e$,
$\lambda_e$ is positive.  The existence of a symplectic structure 
implies that for every $\xi \in \P$, the orientation of
$\Gamma$ associated with $\xi$ has the no-cycle property.
(Proof:  It follows from (\ref{eq:2.8}) that the $\xi$-component of 
$c$ is an $\RR$-valued function on $V_{\Gamma}$ which is
positively oriented with respect to $\xi$.)

\item  \emph{Thom classes:} \label{item:4}
 Fix a vertex $p$ and let $e_1,
  \ldots , e_d$ be the edges containing $p$.  Let $\tau :
  V_{\Gamma} \to S^d(\fg^*)$ be the map which is zero at $q \neq
  p$ and at $p$ is equal to $\alpha_{p, e_1} \ldots \alpha
  _{p,e_d}$.  Then $\tau$ is in 
 $H^{2d}(\Gamma,\alpha)$.

\item \label{item:5} \emph{Sub-objects:} Let $\Gamma_1$ be a 
sub-graph of $\Gamma$ which is compatible with $\alpha$. Then the
inclusion map $j : V_{\Gamma_1} \to V_{\Gamma}$ induces a map
  \begin{displaymath}
    j^* :H(\Gamma, \alpha) \to H(\Gamma_1, \alpha_1), \quad
    \alpha_1 =i^* \alpha \, .
  \end{displaymath}

\item  \emph{Gysin maps:}  \label{item:6}
Suppose that $\Gamma_1$  is compatible 
  with $\alpha$.  The \emph{Thom class} of $\Gamma_1$ is the map 
  $\tau : V_{\Gamma} \to S^{d-r} (\fg^*)$ which is zero on the vertices 
  of $\Gamma$ which are not vertices of $\Gamma_1$ and on
  vertices, $p$, of $\Gamma_1$, is equal to 
  \begin{displaymath}
    \tau (p) = \alpha_{p,e_1} \ldots \alpha_{p,e_s}
  \end{displaymath}
where $s =d-r$ and $e_1 , \ldots , e_s$ are the edges of $\Gamma$ 
at $p$ which don't belong to $\Gamma_1$.

From $\tau$ one gets a Gysin map
\begin{displaymath}
  H^{2k}(\Gamma_1 , \alpha_1) \to H^{2(k+s)}(\Gamma, \alpha)
\end{displaymath}
mapping $f$ to $\tau f$.  (Since $\tau$ is supported on
$V_{\Gamma}$, this map is well-defined.)

\item \emph{The cohomology of blow-ups:} Let $\Gamma$ be a $d$-valent 
graph and $\alpha : I_{\Gamma} \to \fg^*$ an axial function. Let
$p_0$ be a vertex of $\Gamma$ and $(\Gamma^{\#}, \alpha^{\#})$ the 
blow-up of $(\Gamma,\alpha)$ at $p_0$ (See \S \ref{sec:2.1}).
From the blowing down map (\ref{eq:2.4}) one gets a pull-back
map on cohomology
$$\beta^* : H(\Gamma, \alpha) \to H(\Gamma^{\#}, \alpha^{\#})$$
which embeds $H(\Gamma,\alpha)$ as a subring of 
$H(\Gamma^{\#}, \alpha^{\#})$. Moreover the singular locus, $\Gamma_0$,
of $\beta$, is a sub-object of $\Gamma^{\#}$ (in the sense of example
\ref{ex:2} above) and its Thom class,
$$\tau \in H^2 (\Gamma^{\#}, \alpha^{\#}) $$
generates $H(\Gamma^{\#}, \alpha^{\#})$ over the sub-ring
$H(\Gamma, \alpha)$, subject to the relation
$$\tau^d-c_1 \tau^{d-1}+ c_2\tau^{d-2} - ... \pm c_d,$$
the $c_i$'s being the Chern classes of $\Gamma$ for $i <d$ and
$c_d$ being the Thom class of $p_0$.

\end{enumerate}

\subsection{The Atiyah-Bott-Berline-Vergne localization theorem}
\label{sec:2.4new}

  We have just discussed functoriality 
  for sub-objects of $\Gamma$.  What about quotient objects?  To
  take the most extreme case let ``$pt$'' be the trivial
  zero-valent graph consisting of one vertex, $pt$, and no edges
  and let 
  \begin{displaymath}
    \pi : V_{\Gamma} \to pt
  \end{displaymath}
be the constant map.  We have already seen (see \S \ref{sec:2.2} 
item~\ref{item:1}) that
there is a functorial map
\begin{displaymath}
  \pi^* : H^{2k}(pt) \to H^{2k}(\Gamma, \alpha) \, .
\end{displaymath}
However, does there exist a Gysin map
\begin{displaymath}
  \pi_* : H^{2k}(\Gamma, \alpha) \to H^{2(k-d)}(pt) \, ?
\end{displaymath}
Such a map, if it existed, would have to have the following property. 
Let $p$ be a vertex of $\Gamma$ and let $j_p : pt \to V_{\Gamma}$ be 
the map $pt \to p$. Then, by functoriality, $\pi_*$ would have to
satisfy
$$\pi_* (j_p)_* = \mbox{ identity }$$
and, by items \ref{item:4} and \ref{item:6} of \S \ref{sec:2.2}, $\pi_*$
would have to have the form 
\begin{equation}
\label{eq:2.9}
\pi_* f = \sum \frac{f(p)}{\prod \alpha_{p,e}}.
\end{equation}
However, it is by no means obvious that this map is well defined, i.e. 
that the right hand side of (\ref{eq:2.9}) is in $S(\fg^*)$. We will
prove that it is :

\begin{theorem}\label{th:2.1}
$\pi_*$ maps $H^{2k}(\Gamma, \alpha)$ into $S^{k-d}(\fg^*)$.
\end{theorem}

\begin{proof}
Let $f \in H^{2k}(\Gamma, \alpha)$; then $\pi_* f$ can be written as
\begin{equation}
   \label{eq:2.10}
   \pi_* f = \frac{g}{\prod_{j=1}^N \alpha_j}
\end{equation}
where $g \in S^{k-d+N}(\fg^*)$ and $\alpha_1, \cdots , \alpha_N$ are
pair-wise linearly independent. We will show that $\alpha_1$ divides 
$g$.

The vertices of $\Gamma$ can be divided into two categories:

\begin{enumerate}
\item The first subset, $V_1$, contains the vertices 
$p \in V_{\Gamma}$ for which none of the $\alpha_{p,e}$'s is a
multiple of $\alpha_1$

\item The second subset, $V_2$, contains the vertices 
$p \in V_{\Gamma}$ for which there exists an edge $e$ such that
$\alpha_{p,e}$ is a multiple of $\alpha_1$. ( Notice that 
\ref{eq:1.17} implies that there will be exactly one such edge. )

\end{enumerate}

The part of (\ref{eq:2.9}) corresponding to vertices in the first 
category will then be of the form
\begin{equation}
        \label{eq:2.11}
    \sum_{p \in V_1} \frac{f(p)}{\prod \alpha_{p,e}} =
    \frac{g_1}{\prod_{j=2}^N \alpha_j }
\end{equation}
with $g_1 \in S(\fg^*)$.

If $p \in V_2$ then there exists an edge $e$ from $p$ such that 
$\alpha_{p,e} = \lambda \alpha_1$ with 
$\lambda \in \CC - \{ 0 \}$; let $q$ be the other vertex of $e$.
Since $\alpha_{q,e} = - \alpha_{p,e}$
it results that $q \in V_2$ as well and thus the vertices in $V_2$
can be paired as above.

Let
$e_i, i=1, \ldots,d$ and $e'_i , i=1, \ldots ,d$ be the edges
containing $p$ and $q$ respectively, with $e_d=e'_d=e$.  Then the
$e_i$'s can be ordered (cf. \ref{eq:1.18}) so that
\begin{equation}
  \label{eq:2.12}
   \alpha_{p,e_i} \equiv  \alpha_{q,e'_i} \pmod{\alpha_1} \; .
\end{equation}
Also (\ref{eq:1.21}) implies that
\begin{equation}
   \label{eq:2.13}
  f(q) \equiv f(p) \pmod{\alpha_1} \; .
\end{equation}
The part of (\ref{eq:2.9}) corresponding to $p$ and $q$ is given by
\begin{equation}
   \label{eq:2.14}
 \frac{f(p)}{\prod_{j=1}^d \alpha_{p,e_j}} +
 \frac{f(q)}{\prod_{j=1}^d \alpha_{q,e'_j}} =
%
%
\frac{ f(p) \alpha_{q,e'_1}.. \alpha_{q,e'_{d-1}}-
f(q) \alpha_{p,e_1}.. \alpha_{p,e_{d-1}}}{ \lambda \alpha_1
\alpha_{p,e_1}...\alpha_{q,e'_{d-1}}}
\end{equation}
But multiplying together the congruences (\ref{eq:2.12}) and 
(\ref{eq:2.13}) we obtain that $\alpha_1$ divides the numerator of
(\ref{eq:2.10}) so that
\begin{equation}
  \label{eq:2.15}
    \frac{f(p)}{\prod_{j=1}^d \alpha_{p,e_j}} +
    \frac{f(q)}{\prod_{j=1}^d \alpha_{q,e'_j}} =
    \frac{g_{p,q}}{\prod_{j=2}^N \alpha_j}
\end{equation}
with $g_{p,q} \in S(\fg^*)$. Therefore
\begin{equation}
  \label{eq:2.16}
    \sum_{p \in V_2} \frac{f(p)}{\prod \alpha_{p,e}} =
    \frac{g_2}{\prod_{j=2}^N \alpha_j }
\end{equation}
with $g_2 \in S(\fg^*)$. Adding (\ref{eq:2.11}) and (\ref{eq:2.16})
we get that
$$
\frac{g}{\prod_{j=1}^N \alpha_j} =
\frac{ g_1 + g_2}{\prod_{j=2}^N \alpha_j}
$$
with $g_1 + g_2 \in S(\fg^*)$, i.e. that $\alpha_1$ divides $g$.
The same argument can be used to show that each $\alpha_j$ divides 
$g$ and therefore $\pi_* f \in S^{k-d}(\fg^*)$, as desired.

\end{proof}

\subsection{The Kirwan map}
\label{sec:2.3}

From now on we will assume that $(\Gamma, \alpha)$ satisfies the
no-cycle condition. 
Let $\xi$ be an element of $\P$ which gives an orientation of
$\Gamma$ without cycles and let $\phi :V_{\Gamma} \to \RR$
be positively oriented with respect to $\xi$; without loss of 
generality we can assume that $\phi$ is injective.   
For $c \in \RR
- \phi (V_{\Gamma})$, we define the $c$-{\em cross section},  
$\Gamma_c$, of
$\Gamma$ to be the set of edges, $e$, of $\Gamma$ with the
property that, for one of the vertices, $p$, of $e$, $\phi(p)>c$
and for the other vertex, $q$, $\phi (q) <c$.  Let $\g^*_{\xi}$ be
the annihilator of $\xi$ in $\g^*$ and let 
$H^{2k}(\Gamma_c, \alpha)$ be the set of all maps
\begin{displaymath}
  f: \Gamma_c \to S^k (\g^*_{\xi}) \, .
\end{displaymath}
The sum
\begin{displaymath}
  H(\Gamma_c, \alpha) = \bigoplus H^{2k}(\Gamma_c, \alpha)
\end{displaymath}
is a graded ring under point-wise multiplication and we will
define a morphism of graded rings
\begin{displaymath}
  \K_c : H(\Gamma,\alpha) \to H(\Gamma_c,\alpha)
\end{displaymath}
as follows:  For $e \in \Gamma_c$ let $p$ and $q$ be the vertices
of $e$. The projection $\g^* \to \g^*_e$ maps $\g^*_{\xi}$
bijectively onto $\g^*_e$ since $\alpha_{p,e}(\xi) \ne 0$, so one
has a composite map
\begin{displaymath}
  \g^* \to \g^*_e \leftrightarrow \g^*_{\xi}
\end{displaymath}
and hence an induced morphism of graded rings:
\begin{displaymath}
  S(\g^*) \to S(\g^*_e) \leftrightarrow S(\g^*_{\xi}) \, .
\end{displaymath}
If $f$ is in $H(\Gamma,\alpha)$, the images of $f_p$ and 
$f_q$ in $S(\g^*_e)$
are the same by (\ref{eq:1.21}) and hence so are their images in
$S(\g^*_{\xi})$.  We define $\K_c(f)(e)$ to be this common image
and call the map $\K_c$ {\it the Kirwan map}.

Next we will define a morphism, $\gamma_c$, of
$S(\g^*_{\xi})$-modules, mapping $H(\Gamma_c,\alpha)$ 
into the quotient field of
$S(\g^*_{\xi})$.
To define $\gamma_c$, let
$e$, as above, be an element of $\Gamma_c$ and let $p$ and $q$
be the vertices of $e$.  We will assume that $\phi (p) >c$ and
$\phi(q)<c$ and, hence, that $\alpha_{q,e} (\xi)>0$.  Let
\begin{displaymath}
  e^+_i \, , \quad i=1, \ldots , d-1
\end{displaymath}
be the other edges of $\Gamma$ (other than $e$) intersecting at
$p$ and
\begin{displaymath}
  e^-_i \, , \quad i=1, \ldots , d-1
\end{displaymath}
be the other edges of $\Gamma$ intersecting at $q$.  By the
compatibility axiom we can assume that $\alpha_{p,e^+_i}$ and
$\alpha_{q,e^-_i}$ have the same image in $\g^*_e$ and, hence,
under the identification, $\g^*_e \leftrightarrow \g^*_{\xi}$, have
the same image in $\g^*_{\xi}$.  This implies that
\begin{displaymath}
  \alpha^{\#}_{i,e} = : \, \alpha_{p,e^+_i}
  - m^+_{i,e} \alpha_{p,e} = \alpha_{q,e^-_i}
  -m^-_{i,e} \alpha_{q,e}
\end{displaymath}
where
\begin{displaymath}
  m^+_{i,e} = \frac{\alpha_{p,e^+_i}(\xi)}{\alpha_{p,e}(\xi)} \qquad
\mbox{ and } \qquad
  m^-_{i,e} = \frac{\alpha_{q,e^-_i}(\xi)}{\alpha_{q,e}(\xi)}
\end{displaymath}
Let
\begin{displaymath}
  m_e= \alpha_{q,e} (\xi) =- \alpha_{p,e} (\xi)
\end{displaymath}
and note that, since $\phi (p)> \phi(q)$, $m_e= \left| m_e
\right| >0$.  

We now define, for $f \in H^{2k}(\Gamma_c,\alpha)$,
\begin{equation}
  \label{eq:2.17}
  \gamma_c f = \sum_{e \in \Gamma_c} \frac{1}{m_e} \,
  \frac{f(e)}{\prod_i \alpha^{\#}_{i,e}}
\end{equation}
and we define $\p_c$ to be the composition
\begin{equation}
  \label{eq:2.18}
  \p_c =  \gamma_c \K_c \, .
\end{equation}

\begin{theorem}
\label{th:2.2}
$\p_c$ is a map of $H^{2k}(\Gamma,\alpha)$ 
into $S^{k-d+1} (\g^*_{\xi})$.
\end{theorem}

We will prove this by obtaining an explicit ``residue formula''
for $\p_c (f)$, $f \in H^{2k}(\Gamma,\alpha)$.  
This residue formula can be
viewed as a kind of combinatorial version of the 
Jeffrey-Kirwan theorem described in \S\ref{sec:1.3}.  It is also
closely related to the localization theorem proved by Jaap
Kalkman and the first author in \cite{GK}, and the residue results which
we will describe in the next section are mostly taken from
\cite{GK}.

\subsection{Residues}
\label{sec:2.4}

Let $\alpha_1 , \ldots \alpha_d$ be elements of $\g^*$ and 
$\xi$ be 
an element of $\g$ with the property that $\alpha_i (\xi) \ne
0$ for all $i$.  Given $f \in S^k(\g^*)$ we define an element
\begin{equation}
  \label{eq:2.19}
  \Res_{\xi} \frac{f}{\prod \alpha_i}
\end{equation}
of $S^{k-d+1}(\g^*_{\xi})$ as follows: Choose a basis $x, y_1,
.., y_{n-1}$ of $\g^*$ such that $y_1 , .., y_{n-1}$ is a
basis of $\g^*_{\xi}$ and $x(\xi) =1$.  Let
\begin{equation}
  \label{eq:2.20}
  \alpha_i = m_ix - \sum^{n-1}_{j=1} a_{ij} y_j
\end{equation}
and let
\begin{displaymath}
  f(x) = \sum^k_{r=0} x^r f_r (y)
\end{displaymath}
and
\begin{equation}
  \label{eq:2.21}
\beta_i = m^{-1}_i \sum a_{ij} \, y_j \, ;
\end{equation}
then
\begin{displaymath}
  \frac{f}{\alpha_1 \ldots \alpha_d} = (\prod m_i)^{-1}
  x^{-d} \sum x^r f_r \left(1- \frac{\beta_1}{x}\right)^{-1}
  \cdots \left(1- \frac{\beta_d}{x}\right)^{-1} \, .
\end{displaymath}
Now replace $(1- \frac{\beta_i}{x})^{-1}$ by the power series
\begin{displaymath}
  \sum^{\infty}_{k=0} x^{-k} \beta^k_i
\end{displaymath}
and define (\ref{eq:2.19}) to be the coefficient of $x^{-1}$ in
the product on the right.  It is easy to see that this definition
doesn't depend on the choice of $x, y_1, \ldots , y_n$.
If the $\alpha_i$'s are pairwise linearly independent (i.e.,~if,
for $i \ne j$, $\alpha_i$ and $\alpha_j$ are not multiples of
each other) there is a relatively simple formula for
(\ref{eq:2.19}).

\begin{lemma}
  \label{lem:1}
Let $A$ be a graded commutative algebra over $\CC$ and let
$f=f(x)$ be a polynomial in $x$ with coefficients in $A$.  Then
for indeterminants $z_1 , \ldots , z_d$
\begin{equation}
  \label{eq:2.22}
  \Res_x \frac{f(x)}{(x-z_1) \cdots (x-z_d)} =
  \sum^d_{i=1} \frac{f(z_i)}{\prod_{j \ne i} (z_i-z_j)}
\end{equation}
\end{lemma}

\begin{proof}

The conclusion follows immediately from the decomposition in simple
fractions

\begin{displaymath}
\frac{f(x)}{(x-z_1) \cdots (x-z_d)} = F(x)+
\sum^d_{i=1} \frac{f(z_i)}{\prod_{j \ne i} (z_i-z_j)}\frac{1}{x-z_i}
\end{displaymath}
where $F(x)$ is a polynomial term in $x$.

\end{proof}

Let 
$$h=\frac{f}{\prod_{j=1}^{d} (x-z_j)} \quad \mbox{and} \quad
h_j=\frac{f(z_j)}{\prod_{r \neq j} (z_j -z_r)},\; \forall \; j$$

\begin{lemma}\label{lem:2}
$h \in A[x]$ if and only if $\Res_x (x^k h)=0, \; \forall k \geq 0$.
\end{lemma}
 
\begin{proof} From (\ref{eq:2.22}) we get that
$$\Res_x (x^k h)=\sum_{j=1}^d (z_j)^k h_j$$
Then the fact that 
$\Res_x (x^k h)=0, \; \forall k=1,...,d$ can be written as
$$\left( \begin{array}{ccccc}
z_1^1 & \cdots & z_j^1 & \cdots & z_d^1 \\
\vdots & & \vdots & & \vdots \\
z_1^k & \cdots & z_j^k & \cdots & z_d^k \\
\vdots & & \vdots & & \vdots \\
z_1^d & \cdots & z_j^d & \cdots & z_d^d 
\end{array} \right)
\left( \begin{array}{c} h_1 \\ \vdots \\ h_j \\ \vdots \\ h_d
\end{array} \right) = 0$$
Since the corresponding Vandermonde determinant is non-zero we deduce
that $h_1=~\cdots~=~h_d=0$, i.e.
$f(z_j)=0, \; \forall j=1,...,d$, from which we obtain that
$h\in A[x] $; the other implication is clear.                     

\end{proof}  

We will now apply lemma \ref{lem:1} to the evaluation of 
(\ref{eq:2.19}).  Let
\begin{displaymath}
  m_i = \alpha_i (\xi)
\end{displaymath}
and for $i \ne j$ let
\begin{displaymath}
  \alpha^{\#}_{j,i} = \alpha_j - m_{j,i} \alpha_i
\end{displaymath}
with
\begin{displaymath}
  m_{j,i} = \frac{\alpha_j (\xi)}{m_i} \, .
\end{displaymath}
Note that $\alpha^{\#}_{j,i} \in \g^*_{\xi}$ since $\alpha^{\#}_{j,i}
(\xi) = \alpha_j (\xi)- \alpha_j(\xi) =0$.  Let
\begin{displaymath}
  \g^*_i = \g^* / \left\{ c \alpha_i \, ; \, c \in \RR \right\} \, .
\end{displaymath}
The projection map:  $\g^* \to \g^*_i$ is bijective on $\g^*_{\xi}$,
so one gets a composite map
\begin{displaymath}
  \g^* \to \g^*_i \overset{\simeq}{\longrightarrow} \g^*_{\xi}
\end{displaymath}
as in \S\ref{sec:2.3} and hence a ring morphism
\begin{displaymath}
  \K_i : \, S(\g^*) \to S(\g^*_{\xi}) \, .
\end{displaymath}

\begin{theorem}
  \label{th:2.3}
For $f \in S^k(\g^*)$
\begin{equation}
  \label{eq:2.23}
  \Res_{\xi} \frac{f}{\alpha_1 \ldots \alpha_d}
  = \sum_i \frac{1}{m_i}
  \frac{\K_i f}{\prod_{j \ne i} \alpha^{\#}_{j,i}} \, .
\end{equation}

\end{theorem}

\begin{proof}
  With the notations (\ref{eq:2.20})-(\ref{eq:2.21})
  \begin{displaymath}
    \frac{f}{\prod \alpha_i} = (\prod m_k)^{-1}
    \frac{f(x,y)}{\prod (x-\beta_k(y))} \, .
  \end{displaymath}
Thus by lemma \ref{lem:1}
\begin{eqnarray*}
  \Res_{\xi} \frac{f}{\prod \alpha_i}
  &=& (\prod m_k)^{-1} \sum_i
      \frac{f(\beta_i,y)}{\prod_{k \ne i} (\beta_i -
        \beta_k)} \, = \\[2ex]
  &=& \sum_i \frac{1}{m_i}
       \frac{f(\beta_i, y)}{\prod_{k \ne i} m_k (\beta_i - \beta_k)}
\end{eqnarray*}
But $m_k (\beta_i-\beta_k)= \alpha^{\#}_{k,i}$ and the map
\begin{displaymath}
  \K_i : \, S(\g^*) \to S(\g^*_{\xi})
\end{displaymath}
maps $x$ to $\beta_i$ and $y_k$ to itself, so $\K_i f(x,y) =
f(\beta_i, y)$.  Thus the sum on the right is identical with the
right hand side of (\ref{eq:2.22}).

\end{proof}

\subsection{The Jeffrey-Kirwan theorem}
\label{sec:2.5}
We will prove Theorem~\ref{th:2.2} by deducing it from the
following result:

\begin{theorem}
  \label{th:2.4}
For $f \in H^{2k}(\Gamma,\alpha)$
\begin{equation}
  \label{eq:2.24}
  \p_c (f) = \sum_{\phi (p)<c}
  \Res_{\xi} \frac{f_p}{\prod \alpha_{p,e}} \, .
\end{equation}
In particular, $\p_c(f)$ is in $S^{k-d+1}(\g^*_{\xi})$.

\end{theorem}

\begin{proof}
Choose $c =c_0 >c_1>c_2>...>c_N$ in $\RR - \phi (V)$
so that $\phi(V) \subset (c_N, \infty)$ and  $\forall r \geq 0$
there is exactly one vertex $p_r$ with 
$\phi (p_r) \in ( c_{r+1}, c_r)$. 
Inspection of (\ref{eq:2.18}) shows that
  \begin{equation}
    \label{eq:2.25}
    \p_{c_r} (f)- \p_{c_{r+1}}(f) = 
  \sum \frac{1}{m_i}
    \frac{\K_i f}{\prod_{j \ne i} \alpha^{\#}_{j,i}}  \, ,
  \end{equation}
where $e_i$, $i=1, \ldots ,d$, are the edges of $\Gamma$
containing $p_r$ and $\alpha_i = \alpha_{p_r,e_i}$.  On the other
hand, by Theorem~\ref{th:2.3}, the right hand side of
(\ref{eq:2.25}) is just
\begin{displaymath}
  \Res_{\xi} \frac{f(p_r)}{\prod \alpha_{p_r,e}} \, .
\end{displaymath}
The conclusion follows since $\p_{c_N} (f)=0$.

\end{proof}

\begin{corollary}
 \label{cor:2.1}
If $\pi_*$ is the map given by (\ref{eq:2.9}) then, 
 for $f \in H^{2k}(\Gamma,\alpha)$,
\begin{equation}
  \label{eq:2.26}
  \Res_{\xi}(\pi_*f)=0
\end{equation}
\end{corollary}

We conclude this section 
by observing that combining corollary \ref{cor:2.1} and
lemma \ref{lem:2} we obtain a new proof of theorem \ref{th:2.1}
for graphs that satisfy the no-cycle condition:

Let $f \in H^{2k}(\Gamma,\alpha)$. Then, as in (\ref{eq:2.10}),
\begin{equation*}
   \pi_* f = \frac{g}{\prod_{j=1}^N \alpha_j}
\end{equation*}
where $g \in S^{k-d+N}(\fg^*)$ and $\alpha_1, \cdots , \alpha_N$ are
pair-wise linearly independent.

Let $\xi$ generate an orientation of $\Gamma$ with no cycles and choose
$\theta \in \g^*$ such that 
$\theta (\xi)=1$ and $\theta$ is not equal to any 
of $\alpha_1,\cdots,\alpha_N$. 
Then $\theta^r f \in H^{2(k+r)}(\Gamma,\alpha)$ 
and 
\begin{displaymath}
    \pi_* (\theta^r f)= \frac{\theta^r g}{\prod_{i=1}^N \alpha_i}
  \end{displaymath}
But (\ref{eq:2.26}) implies that 
$\Res_{\xi}(\pi_*(\theta^r f))=0 \; \forall r \geq 0$ and it follows 
now from lemma \ref{lem:2} that $\pi_*f \in S^{k-d}(\g^*)$.

\subsection{The Betti numbers of the pair ($\Gamma, \alpha$)}
\label{sec:2.6}

For $\xi \in \P$ and $p \in V_{\Gamma}$, let 
$\sigma_p = \sigma_p(\xi)$ be the number of edges $e$, with one
vertex $p$, for which $\alpha_{p,e} (\xi) <0$. Let $\beta_k$ be the
number of vertices $p$ with $\sigma_p = k$. Since $\sigma_p$ 
depends on $\xi$, it is surprising to find that these ``Betti
numbers'' don't.

\begin{theorem}
\label{th:2.5}
The $\beta_k$'s don't depend on $\xi$; i.e. they are combinatorial 
invariants of $(\Gamma, \alpha)$.
\end{theorem}

\begin{proof}
Let $\P_i, \; i=1,...,N$, be the connected components of $\P$
and consider an $(n-1)$-dimensional wall separating two adjacent 
$\P_i$'s. This wall is defined by an equation of the form
\begin{equation}
\label{eq:2.27}
\alpha_{p,e}(\xi) =0
\end{equation}
for some $(p,e) \in I_{\Gamma}$. Let $q$ be the other vertex of $e$
and lets compute the changes in $\sigma_p$ and $\sigma_q$ as $\xi$
passes through this wall: Let $e_i, \; i=1,...,d$ be the edges 
meeting at $p$ and $e'_i, \; i=1,...,d$ be the edges meeting at $q$
( with $e_d =e'_d=e$ ).  By (\ref{eq:1.18}) we can order the $e_i$'s 
so that, for $i \leq d-1$,
\begin{displaymath}
\alpha_{p,e_i} = \alpha_{q,e'_i} +c_i \alpha_{p,e} \; .
\end{displaymath}
From (\ref{eq:1.17}) follows that for every $i=1,...,d-1$,
\begin{displaymath}
\mbox{dim}\, (\,\mbox{ker}\, \alpha_{p,e} \cap \mbox{ ker}\, 
\alpha_{p,e_i}) =n-2
\end{displaymath}
Therefore there exists $\xi_0$ such that 
$\alpha_{p,e}(\xi_0)=0$ but 
$\alpha_{p,e_i}(\xi_0) = \alpha_{q,e'_i}(\xi_0) 
\neq 0$, for all $i=1,...,d-1$.

Then there exists a neighborhood $U$ of $\xi_0$ in $\fg$ such that
for all $i=1,...,d-1$ and $\xi \in U$,  
$\alpha_{p,e_i}(\xi)$ and  $ \alpha_{q,e'_i}(\xi)$
have the same sign
and this common sign doesn't depend on $\xi \in U$.
Such a neighborhood will intersect both regions created by the 
wall (\ref{eq:2.27}). Now suppose that $\xi \in U$ and that
$r$ of the numbers $\alpha_{p,e_i}(\xi), \; i=1,...,d-1,$ are 
negative. Since $\alpha_{p,e}(\xi) = - \alpha_{q,e}(\xi)$, it
follows that for $\alpha_{p,e}(\xi)$ positive
\begin{displaymath}
\sigma_p =r \quad \mbox{ and } \quad \sigma_q =r+1
\end{displaymath}
and for $\alpha_{p,e}(\xi)$ negative 
\begin{displaymath}
\sigma_p =r+1 \quad \mbox{ and } \quad \sigma_q =r
\end{displaymath}
In either case, as $\xi$ passes through the wall (\ref{eq:2.27}),
the Betti numbers don't change. 

\end{proof}
(For this simple and beautiful
proof of the well-definedness of the Betti numbers we are indebted
to Ethan Bolker.)

\subsection{Betti numbers and cohomology}
\label{sec:2.9}

Simple examples show that the formula (\ref{eq:1.35}) is not true 
for an arbitrary graph-axial function pair $(\Gamma, \alpha)$.
However, we will prove that if $(\Gamma, \alpha)$ has the no-cycle 
property for some $\xi \in \P$ then the equality (\ref{eq:1.35})
can be replaced by the inequality
\begin{equation}
\label{eq:2.28}
\mbox{ dim }H^{2k}(\Gamma, \alpha) \leq 
\sum \beta_r \mbox{ dim } S^{k-r}(\fg^*)
\end{equation}

In addition we will show that, for $k$ large,
\begin{equation}
\label{eq:2.29}
\mbox{ dim } H^{2k}(\Gamma, \alpha) =
\sum \beta_r \mbox{ dim } S^{k-r}(\fg^*)+ O(k^{n-3}).
\end{equation}
(Note that since
$$\mbox{ dim } S^k(\fg^*) = 
\left( \begin{array}{c} k+n-1 \\ n-1 \end{array} \right) =
\frac{1}{(n-1)!} \left( k^{n-1}+
\left( \begin{array}{c} n\\2 \end{array} \right) k^{n-2} +
O(k^{n-3}) \right)$$
the first term on the right hand side is strictly greater 
than the error term.) In particular, if $n=2$, the formula 
(\ref{eq:2.29}) asserts that
\begin{equation}
\label{eq:2.30}
\mbox{ dim } H^{2k}(\Gamma, \alpha) =
\sum \beta_r \mbox{ dim } S^{k-r}(\fg^*)
\end{equation}
for all $k$ greater than some fixed $k_0$.

\begin{proof}
Let $\alpha_i \in \g^*,\; i=1,...,N$, be a pairwise linearly 
independent set of vectors with the property that every one 
of the vectors $\alpha_{p,e}$, $(p,e) \in I_{\Gamma}$, is a 
multiple of a vector in this set. Let $I$ be the graded ideal
in $S(\g^*)$ generated by the monomials 
$$g_i=\alpha_1 \cdots \widehat{\alpha_i} \cdots \alpha_N.$$

\begin{lemma} 
\label{lem:3}
The algebraic dimension of the quotient ring $S(\g^*)/I$ 
is $n-2$.
\end{lemma}

\begin{proof} This follows trivially from the fact that the 
algebraic variety defined by $I$ is the union of the sets
$\alpha_i=\alpha_j =0,\; i \neq j$.
\end{proof}

As a corollary of this lemma we get the bound
\begin{equation}
\label{eq:2.31}
\mbox{ dim } S^k(\fg^*)/I^k = O(k^{n-3}).
\end{equation}
Now let $\phi : V_{\Gamma} \to \RR$ be a strictly monotone function
which is positively oriented with respect to $\xi$ and let 
$H_c(\Gamma,\alpha)$ be the subring of $H(\Gamma,\alpha)$ consisting
of all maps $f : V_{\Gamma} \to S(\g^*)$ with support on the set 
$\phi \geq c$. Let $p \in V_{\Gamma}$ with $\phi (p) =c$ and suppose 
that there are no points $q \in V_{\Gamma}$ with $\phi (q)$ on
the interval $(c,c')$. Letting $\sigma_p=r$ we will prove the Morse
inequality 
\begin{equation}
\label{eq:2.32}
\mbox{ dim } H_c^{2k}(\Gamma,\alpha) -
\mbox{ dim } H_{c'}^{2k}(\Gamma,\alpha) 
\leq \mbox{ dim }S^{k-r}(\fg^*)
\end{equation}
and an inequality in the opposite direction:
\begin{equation}
\label{eq:2.33}
\mbox{ dim } H_c^{2k}(\Gamma,\alpha) -
\mbox{ dim } H_{c'}^{2k}(\Gamma,\alpha) 
\geq \mbox{ dim }I^{k-r}
\end{equation}
To prove (\ref{eq:2.32}) let $e_i,\; i=1,...,d$, be the edges of $\Gamma$
containing $p$ and let $\alpha_i=\alpha_{p,e_i}$. We will order
the $\alpha_i$'s so that $\alpha_i(\xi)<0$ for $1 \leq i \leq r$ and
$\alpha_i(\xi) >0$ for $r+1 \leq i \leq d$. Then if 
$f \in H_c^{2k}(\Gamma,\alpha)$, $f(p)$ must be a multiple of 
$\alpha_1,..., \alpha_r$ so the image of the restriction map
\begin{equation}
\label{eq:2.34}
H_c^{2k}(\Gamma,\alpha) \to S^k(\fg^*), \quad f \to f(p),
\end{equation}
is contained in $S^{k-r}(\fg^*)\alpha_1\cdots \alpha_r$. 
Since the kernel
of this map is $H_{c'}^{2k}(\Gamma,\alpha)$, this proves 
(\ref{eq:2.32}). We will prove the inequality (\ref{eq:2.33}) by 
showing that if $h \in I^{k-r}$ then $h \alpha_1 \cdots \alpha_r$ 
is in the image of (\ref{eq:2.34}). Indeed, if $h \in I^{k-r}$ then
$h \alpha_1 \cdots \alpha_r$ can be written as a sum
$$\sum_{i=r+1}^N h_i \alpha_1 
\cdots \widehat{\alpha_i} \cdots \alpha_N. $$
Let $p_j$ be the vertex joined to $p$ by $e_j$ for $j=r+1,\cdots,d$, 
and, for fixed $j_0 \in \{ r+1,..., d \}$, 
define $f: V_{\Gamma} \to S^k(\g^*)$ to be the map which takes the value 
$h \alpha_1 \cdots \alpha_r$ at $p$, the value 
$h_{j_0} \alpha_1 \cdots \widehat{\alpha_{j_0}} 
\cdots \alpha_N$ at $p_{j_0}$ and zero
elsewhere. It is easily checked that $f \in H_c^{2k}(\Gamma,\alpha)$ and
$f(p)=h \alpha_1 \cdots \alpha_r$. This proves (\ref{eq:2.33}).

Next let $c$ and $c'$ be any pair of real numbers with $c < c'$.
From (\ref{eq:2.32})-(\ref{eq:2.33}) one gets, by a simple induction,
the Morse inequalities 
$$\mbox{ dim } H_c^{2k}(\Gamma,\alpha) -
\mbox{ dim } H_{c'}^{2k}(\Gamma,\alpha) 
\leq \sum \beta_r(c,c') \mbox{ dim } S^{k-r}(\fg^*)$$
and
$$\mbox{ dim } H_c^{2k}(\Gamma,\alpha) -
\mbox{ dim } H_{c'}^{2k}(\Gamma,\alpha) 
\geq \sum \beta_r(c,c') \mbox{ dim } I^{k-r}$$
where $\beta_r(c,c')$ is the number of points $p \in V_{\Gamma}$ 
with $\sigma_p=r$ and $c \leq \phi(p) < c'$. In particular, for
$c'>>0$ and $c<<0$, one gets from these estimates and from 
(\ref{eq:2.31}) the inequalities (\ref{eq:2.28}) and (\ref{eq:2.29}).

\end{proof}

\subsection{The role of the zeroth Betti number}
\label{sec:2.10}

An example of a graph-axial function pair that fails to 
satisfy (\ref{eq:1.35}) is the $d=n=2$ example described at the end of
\S \ref{sec:2.1}. The graph in this example is a connected graph; 
so its topological zeroth Betti number is 1. However, its graph
theoretical zeroth Betti number, $\beta_0$, is $N$. A simple computation
shows that, for this example, the identity (\ref{eq:1.35}) holds for all 
$k >0$. But for $k=0$ the left hand side of (\ref{eq:1.35}) is 1 
(since the graph is connected) whereas the right hand side, $\beta_0$, 
is $N$. From this example one can generate examples of graph-axial 
function pairs, $(\Gamma,\alpha)$, for which the estimate 
(\ref{eq:2.29}) is ``best possible'' by taking Cartesian products of 
this graph with graphs which do satisfy (\ref{eq:1.35}).

However, by making some additional assumptions on the pair 
$(\Gamma,\alpha)$ one can considerably improve (\ref{eq:2.29}). 
The assumptions we will make are of two kinds:
\begin{enumerate}
\item To avoid the problem posed by the example we have just described, 
we will assume that the zeroth Betti numbers of certain connected
subgraphs of $\Gamma$ are equal to 1.
\item For every $p \in V_{\Gamma}$ we will make certain 
``general position'' hypothesis about the vectors 
$\alpha_{p,e}$, $(p,e) \in I_{\Gamma}$. To formulate these hypothesis
we introduce the following refinement of the notion of ``pairwise 
linearly independent'':
\end{enumerate}

\begin{definition}
A collection of vectors $\alpha_i \in \fg^*$, $i=1,...,N$, is 
$l$-{\em independent} if, for every sequence 
$1 \leq i_1 < i_2 < ... < i_l \leq N$, the vectors 
$\alpha_{i_1},...,\alpha_{i_l}$ are linearly independent.
\end{definition}

Now let $(\Gamma,\alpha)$ be a graph-axial function pair
which satisfies the no-cycle condition for some $\xi \in \P$. The 
main result of this section is the following sharpening of
(\ref{eq:2.29}):

\begin{theorem}
\label{th:2.6}
Suppose the following hypotheses hold:
\begin{enumerate}
\item For every subspace $\fh$ of $\fg$ of codimension strictly 
less than $l$, the zeroth Betti number of the connected components of
$\Gamma_{\fh}$ are equal to 1.
\item For every vertex $p$ of $\Gamma$, the vectors $\alpha_{p,e}$,
$(p,e) \in I_{\Gamma}$, are $l$-independent.
\end{enumerate}
Then:
\begin{equation}
\label{eq:2.35}
\mbox{ dim } H^{2k}(\Gamma,\alpha) = 
\sum \beta_r \mbox{ dim } S^{k-r}(\fg^*) + O(k^{n-1-l}).
\end{equation}
\end{theorem}

\begin{remark*}
For $l=2$ the above conditions are always satisfied; 
(\ref{eq:2.29}) is the particular case of (\ref{eq:2.35})
corresponding to $l=2$.
\end{remark*}

For $l=n$ this theorem says that the left hand side of (\ref{eq:2.35})
is equal to the first term on the right for $k$ greater than some fixed 
$k_0$. This result can be slightly improved.

\begin{theorem}
\label{th:2.7}
If the hypotheses of theorem \ref{th:2.6} hold for $l=n$ then
\begin{equation}
\label{eq:2.36}
\mbox{ dim } H^{2k}(\Gamma,\alpha) = 
\sum \beta_r \mbox{ dim } S^{k-r}(\fg^*)
\end{equation}
for $k > d-n$.
\end{theorem}

We will prove these two results by refining the Morse inequalities 
(\ref{eq:2.33}). For this we will need the following generalization of 
Lemma \ref{lem:3} of \S \ref{sec:2.9}:

\begin{lemma}
\label{lem:4}
Let $\gamma_1,...,\gamma_N$ be a collection of vectors in $\fg^*$ which
are $l$-independent and let $I_l$ be the ideal in $S(\fg^*)$ generated
by the monomials
\begin{equation}
\label{eq:2.37}
\frac{\gamma_1 \cdots \gamma_N}{\gamma_{i_1} \cdots \gamma_{i_{l-1}}},
\qquad 1 \leq i_1 < \cdots i_{l-1} \leq N.
\end{equation}
Then the algebraic dimension of $S(\fg^*)/I_l$ is 
$n-l$. Moreover, if $N \geq n$ and $n=l$
then $S^m(\fg^*) = I_l^m$ for $m > N-n$.
\end{lemma}

\noindent 
(For the proof of this lemma see the appendix at the end of this section.)

\begin{proof}{\it (of theorem \ref{th:2.6})}
Let $\phi : V_{\Gamma} \to \RR$ be a strictly monotone function 
which is positively oriented with respect to $\xi$. Let $p$ be an 
arbitrary vertex of $\Gamma$, let $c=\phi(p)$ and assume that there are
no points $q \in V_{\Gamma}$ with $c < \phi (q) < c'$.
Let $e_1,...,e_d$ be the edges of $\Gamma$ containing $p$ and let 
$\alpha_i= \alpha_{p,e_i}$. We can order these vectors so that 
$\alpha_i(\xi) <0$ for $1 \leq i \leq r$ and  
$\alpha_i(\xi) >0$ for $r+1 \leq i \leq d$, 
where $r=\sigma_p$. Let $N=d-r$.

Suppose $l \leq N$. We will prove that
\begin{equation}
\label{eq:2.38}
\mbox{ dim } H_c^{2k}(\Gamma,\alpha) -
\mbox{ dim } H_{c'}^{2k}(\Gamma,\alpha) 
\geq \mbox{ dim }I_l^{k-r}\; ,
\end{equation}
where $I_l$ is the ideal of $S(\fg^*)$ constructed as in lemma
\ref{lem:4} using the $N=d-r$ vectors 
$\alpha_{r+1},...,\alpha_d$ 
(which are $l$-independent, by hypothesis 2), i.e.
$I_l$ is the ideal in $S(\fg^*)$ generated by the monomials
$$\frac{\alpha_{r+1} \cdots \alpha_d}
{\alpha_{i_1}\cdots \alpha_{i_{l-1}}} \; , \quad 
r+1 \leq i_1 < ... < i_{l-1} \leq N.$$
To show (\ref{eq:2.38}), consider (as in \S \ref{sec:2.9}) 
the restriction map
\begin{equation}
\label{eq:2.39}
H_c^{2k}(\Gamma,\alpha) \longrightarrow 
S^{k-r}(\fg^*)\alpha_1 \cdots \alpha_r
\end{equation}
given by $f \to f(p)$. The kernel of this map is
$H_{c'}^{2k}(\Gamma,\alpha)$; so, to prove (\ref{eq:2.38}) it suffices to 
show that the image of this map contains 
$I_l^{k-r}\alpha_1 \cdots \alpha_r$. Consider the set of vectors
$ \alpha_{i_1}, ..., \alpha_{i_{l-1}}, \; r+1 \leq i_1<...<i_{l-1} \leq d$.
These vectors are linearly independent; so their annihilator, $\fh$, 
is of codimension $l-1$ in $\fg$. Let $\Gamma_1$ be the connected
component of $\Gamma_{\fh}$ containing $p$. Since the numbers 
$\alpha_{i_1}(\xi), ..., \alpha_{i_{l-1}}(\xi)$ are greater than zero,
$p$ is a local minimum point for the restriction of $\phi$ to $\Gamma_{\fh}$; 
so, by hypothesis 1, $p$ is also a global minimum. Therefore 
the vertices of $\Gamma_1$ are contained in the set
$\phi \geq c$, and, hence the Thom class, $\tau_1$, of $\Gamma_1$, is 
supported in this set, i.e. is an element of $H_c(\Gamma,\alpha)$.
However, at $p$, $\tau_1 $ is equal to 
$$\left( \frac{\alpha_{r+1} \cdots \alpha_d}
{\alpha_{i_1} \cdots \alpha_{i_{l-1}}} \right) 
\alpha_1 \cdots \alpha_r.$$
This argument shows that, for all generators, $f$, of $I_l$, 
$f\alpha_1\cdots \alpha_r$ is in the image of the restriction
map (\ref{eq:2.39}). Hence the image of this restriction map contains
$I_l^{k-r}\alpha_1 \cdots \alpha_r$, which proves (\ref{eq:2.38}).

Now suppose that $l > N$. In this case we can simply take $\fh$ to be
the annihilator of $\alpha_{r+1},...,\alpha_d$, and, 
by the same argument as
above, conclude that the image of the restriction map is equal to 
$S^{k-r}(\fg^*)\alpha_1 \cdots \alpha_r$. Hence for $l > N$,
\begin{equation}
\label{eq:2.40}
\mbox{ dim } H_c^{2k}(\Gamma,\alpha) -
\mbox{ dim } H_{c'}^{2k}(\Gamma,\alpha) 
= \mbox{ dim }S^{k-r}(\fg^*).
\end{equation}
The proof of the estimate (\ref{eq:2.35}) via (\ref{eq:2.38}),
(\ref{eq:2.40}) and lemma \ref{lem:4} is the same as the proof of the 
estimate (\ref{eq:2.29}) via (\ref{eq:2.33}) and lemma \ref{lem:3}.
We will omit details.

\end{proof}

\begin{proof}({\it of theorem \ref{th:2.7}})
Let $l=n$. If $N<n$ the equality (\ref{eq:2.40}) holds for all $k$
(as we have just seen). If $N \geq n$ then, by lemma \ref{lem:4} and
by (\ref{eq:2.38}), the equality (\ref{eq:2.40}) holds if
$k-r > N-n$, i.e. if $k > d-n$. Thus for $l=n$, (\ref{eq:2.35}) can be
sharpened to (\ref{eq:2.36}).

\end{proof}

\subsection*{Appendix: The proof of Lemma \ref{lem:4}}
\label{appendix}

We will prove by induction that the algebraic variety defined by $I_l$
is the union of the sets
$$\gamma_{i_1}= \cdots = \gamma_{i_l}=0, \quad 
1 \leq i_1< \cdots < i_l \leq N.$$
Let $x$ be a point on this variety. Since $I_{l-1} \subset I_l$, the
variety defined by $I_{l-1}$ contains the variety defined by $I_l$; so
it follows by induction that there exists a sequence
$1 \leq i_1 < \cdots < i_{l-1} \leq N$ with 
$\gamma_{i_1}= \cdots = \gamma_{i_{l-1}}=0$ at $x$. However, since
$I_l$ contains the quotient of $\gamma_1 \cdots \gamma_N$ by
$\gamma_{i_1} \cdots \gamma_{i_{l-1}}$, there exists some 
$j \neq i_1,...,i_{l-1}$ such that $\gamma_j=0$ at $x$. This proves
the first assertion of lemma \ref{lem:4}.
Now let $l=n$. We will prove that
$$S^m(\fg^*)=I_n^m$$
for $m>N-n$ by a double induction on $n$ and $N$.

The equality above is true if $N=n$ or if $n=1$, as can be easily 
checked. Consider now a pair $(N,n)$ with $N > n$.
We now assume that the assertion is true for $(N-1,n-1)$ and for 
$(N-1,n)$.
 
Let $\fh$ be the annihilator of $\gamma_N$ in $\fg$. The restriction 
map
$$S(\fg^*) \longrightarrow S(\fh^*) $$
maps $\gamma_1,...,\gamma_{N-1}$ onto vectors $\beta_1,...,\beta_{N-1}$,
which are $(n-1)$-independent in $\fh^*$; so, by induction, every element 
of $S^m(\fh^*)$, $m > N-n$, is in the ideal generated by the monomials
$$\frac{\beta_1 \cdots \beta_{N-1}}{\beta_{i_1} \cdots \beta_{i_{n-2}}},
\quad 1 \leq i_1 < ... < i_{n-2} \leq N-1.$$
Since the kernel of the map $S(\fg^*) \to S(\fh^*)$ is the ideal
generated by $\gamma_N$, it follows that for $m > N-n$, every element of
$S^m(\fg^*)$ can be written as a linear combination of
$$\frac{\gamma_1 \cdots \gamma_{N-1}\gamma_N}
{\gamma_{i_1} \cdots \gamma_{i_{n-2}}\gamma_N}, \quad
1 \leq i_1 < ...< i_{n-2} \leq N-1,$$
with polynomial coefficients, plus a term of the form $f\gamma_N$,
$f \in S^{m-1}(\fg^*)$. By induction the
theorem is true in dimension $n$ for the vectors 
$\gamma_1,...,\gamma_{N-1}$. Then $f$ is in the  ideal generated by the 
monomials
$$\frac{\gamma_1 \cdots \gamma_{N-1}}
{\gamma_{i_1} \cdots \gamma_{i_{n-1}}}, \quad
1 \leq i_1 < ... < i_{n-1} \leq N-1$$
and hence $f\gamma_N$ is in the ideal generated by
$$\frac{\gamma_1 \cdots \gamma_N}
{\gamma_{i_1} \cdots \gamma_{i_{n-1}}}, \quad
1 \leq i_1 < ... < i_{n-1} \leq N-1.$$

\end{document}